\documentclass[reqno]
{amsart}

\usepackage[T2A]{fontenc}
\usepackage[cp1251]{inputenc}
\usepackage[english,russian]{babel}
\usepackage{amsmath}
\usepackage{amsfonts}
\usepackage{amssymb}
\usepackage{graphicx}
\usepackage{xcolor}
\usepackage{mathrsfs}

\usepackage{enumerate}

\usepackage{upgreek}

\usepackage[pdftex,unicode,colorlinks,linkcolor=red,
citecolor=red,bookmarksopen,pdfhighlight=/N]{hyperref}

\numberwithin{equation}{subsection}

\usepackage[hyper]{amsbib}

\usepackage[hyper]{amsbib}

\theoremstyle{plain} 
\newtheorem{theorem}{Теорема} 
\newtheorem*{maintheorem}{Основная теорема}

\newtheorem{theoB}{Теорема B\hspace{-5pt}} 

\newtheorem*{thMR}{Теорема Мальявена\,--\,Рубела}
\newtheorem*{theoKh}{Теорема Kh} 

\newtheorem*{thWHL}{Теорема Вейерштрасса\,--\,Адамара\,--\,Линделёфа}
\newtheorem*{thWHLB}{Теорема Вейерштрасса\,--\,Адамара\,--\,Линделёфа\,--\,Брело}

\newtheorem{lemma}{Лемма}[section]
\newtheorem{propos}{Предложение}[section] 

\theoremstyle{definition}
\newtheorem{definition}{Определение}
\newtheorem{remark}{Замечание}[section] 
\newtheorem{example}{Пример}[section]

\theoremstyle{plain} 
\newtoks\thehProclaim 
\newtheorem*{Proclaim}{\the\thehProclaim}

\theoremstyle{definition} 
\newtoks{\thehRemark} \newtheorem*{Remark}{\the\thehRemark}

\renewcommand{\leq}{\leqslant} 
\renewcommand{\geq}{\geqslant}
\newcommand{\rad}{\text{\tiny\rm rad}}
\newcommand{\bal}{\rm {bal}}
\newcommand{\Bal}{\rm {Bal}}
\newcommand{\RR}{\mathbb{R}} 
\newcommand{\CC}{\mathbb{C}} 
\newcommand{\NN}{\mathbb{N}} 

\newcommand{\DD}{\mathbb{D}}

\DeclareMathOperator{\har}{har} 
\DeclareMathOperator{\Hol}{Hol} 

\DeclareMathOperator{\Meas}{{\mathcal M}} 
\DeclareMathOperator{\Zero}{Zero} 
\DeclareMathOperator{\sbh}{sbh} 
\DeclareMathOperator{\dsbh}{\text{$\updelta${\rm -sbh}}} 
\DeclareMathOperator{\supp}{supp} 
\DeclareMathOperator{\type}{type} 
\DeclareMathOperator{\ord}{ord}
\DeclareMathOperator{\up}{\text{\rm \tiny up}}
 \DeclareMathOperator{\lw}{\text{\rm \tiny low}}
\DeclareMathOperator{\rght}{\text{\rm \tiny rh}}
 \DeclareMathOperator{\lft}{\text{\rm \tiny lh}}
\DeclareMathOperator{\dd}{\,{\mathrm d\!}}
 
\renewcommand{\Re}{{\rm Re \,}}
\renewcommand{\Im}{{\rm Im \,}}

\begin{document} 
\title[Обобщение и развитие теоремы Мальявена\,--\,Рубела \dots]{Обобщение и развитие теоремы Мальявена\,--\,Рубела о малости роста целых функций экспоненциального типа с заданными нулями}

\author{Б.\,Н.~Хабибуллин, А.\,Е.~Салимова} 
	
\address{факультет математики и ИТ\\ Башкирский государственный университет\\ 450074, г. Уфа\\ ул. Заки Валиди, 32\\ Башкортостан\\
Россия}
	
	\email{Khabib-Bulat@mail.ru}
\email{anegorova94@bk.ru}
	
	
\begin{abstract} 

Ранее нами была  развита и разработана техника выметания соответственно  рода $0$ и рода  $q\in \mathbb N:=\{1,2,\dots\}$ меры или заряда и ($\updelta$-)су\-б\-г\-а\-р\-м\-о\-н\-и\-ч\-е\-с\-к\-ой функции конечного порядка на произвольную замкнутую систему лучей $S$ с вершиной в нуле на комплексной плоскости $\mathbb C$. В данной статье  мы используем только два рода выметания мер и зарядов, а также  субгармонических функций конечного типа при порядке $1$ и их разностей. Первый   --- это классическое выметание рода $q=0$ на систему из четырёх замкнутых лучей: положительная и отрицательная вещественные и мнимые полуоси 
$\mathbb R^+$, $-\mathbb R^+$,  $i\mathbb R^+$, $-i\mathbb R$. Второй   
--- двустороннее выметание рода $q=1$ из открытых правой и левой полуплоскостей $\mathbb C_{\text{\rm \tiny rh}}$ и $\mathbb C_{\text{\rm \tiny lh}}$ на мнимую ось $i\mathbb R$.  
Классическая теорема-критерий   Мальявена\,--\,Рубела даёт законченные условия существования целой функции экспоненциального типа (пишем ц.ф.э.т.) $f\not\equiv 0$, обращающейся в нуль на заданной   \textit{положительной\/} последовательности ${\sf Z}=\{{\sf z}_k\}_{k\in \mathbb N}\subset \mathbb R^+$ и удовлетворяющей ограничению $|f|\leq |g|$ на $i\mathbb R$, где $g$ ---  ц.ф.э.т., обращающаяся в нуль на  \textit{положительной\/} последовательности  ${\sf W}=\{{\sf w}_k\}_{k\in \mathbb N}\subset \mathbb R^+$. 
Сочетание указанных выше специальных выметаний рода $q=0$ и $q=1$
позволяет распространить теорему  Мальявена\,--\,Рубела на произвольные  \textit{комплексные\/} последовательности ${\sf Z}=\{{\sf z}_k\}_{k\in \mathbb N}\subset \mathbb C$, отделённые какой-либо парой вертикальных углов от мнимой оси $i\mathbb R$,  со значительно более общими ограничениями  $\ln |f|\leq M$ на  мнимой оси $i\mathbb R$, где $M$ ---  произвольная субгармоническая функция конечного типа при порядке $1$.
\end{abstract}
	
\thanks{Исследование выполнено за счёт гранта Российского
научного фонда (проект № 18-11-00002).}	
\date{1 мая 2020 г.}
	
\maketitle

\tableofcontents

\section{Введение}\label{s10}

Мы используем обозначения, терминологию, определения и соглашения 
 из  \cite{KhI} и \cite{KhII}.  Так,
\begin{equation*}\label{C+-}
\begin{split}
\CC_{\rght}:= \{z\in \CC \colon \Re z>0\},&\quad 
{\CC_{\lft}}:=\{z\in \CC\colon \Re z<0\},\\
\CC^{\up}:= \{z\in \CC \colon \Im z>0\},&\quad \CC_{\lw}:=\{z\in \CC\colon \Im z<0\}, \\
\CC_{\overline \rght}:=\CC_{\rght}\cup i\RR,
\quad \CC_{\overline \lft}:={\CC}_{\lft}\cup i\RR,&
\quad \CC^{\overline \up}:=\CC^{\up}\cup \RR, 
\quad \CC_{\overline \lw}:=\CC_{\lw}\cup \RR,
\end{split}
\end{equation*} 
--- соответственно {\it открытые  правая и левая, верхняя и нижняя полуплоскости,\/} а  
также их {\it замыкания  в комплексной плоскости\/} $\CC$.

\subsection{Предшествующие результаты}
\textit{\underline{Всюду далее}}
\begin{equation}\label{Z}
{\sf Z}=\{{\sf z}_k\}_{k= 1,2,\dots}\subset \CC, 
\quad 
{\sf W}=\{{\sf w}_k\}_{k:=1,2,\dots}\subset \CC
\end{equation} ---
две последовательности {\it без предельных точек в комплексной плоскости $\CC$}. Последовательности ${\sf W}$ 
сопоставляется её {\it радиальная считающая функция\/}
$n_{\sf W}^{\rad}$, определённая на {\it положительной полуоси\/} $\RR^+:=\{r\in \RR\colon r\geq 0\}$ {\it вещественной оси\/} $\RR$ как 
\begin{equation}\label{nZrad}
 n_{\sf W}^{\rad}\colon t\underset{\text{\rm \tiny $t\in \RR^+$}}{\longmapsto} \sum_{|{\sf w}_k|\leq t}1.
\end{equation}
Для ${\sf W}$ предполагается  конечность \textit{верхней плотности при порядке\/} $1$:
\begin{equation}\label{Zt}
\type_1[{\sf W}] :=\limsup_{t\to \infty} \frac{n_{\sf W}^{\rad}(t)}{t}<+\infty. 
\end{equation}  

\begin{thMR}[{\rm \cite[теорема 4.1]{MR}, см. также \cite[гл.~22, основная теорема]{RC}}]\label{MR1} Пусть 
последовательности ${\sf W}\subset \RR^+$,  ${\sf Z}\subset \RR^+$  положительны  и $\type_1[{\sf Z}] <+\infty$ дополнительно к \eqref{Zt}. 

Тогда  эквивалентны  следующие три утверждения:
\begin{enumerate}[{\rm (i)}]
\item\label{fgi} для любой ц.ф.э.т. $g\neq 0$, обращающейся в нуль на ${\sf W}$, т.е. при $g({\sf W})=0$,  найдётся ц.ф.э.т. $f\neq 0$, для которой
$f({\sf Z})=0$ и 
\begin{equation}\label{fgiR}
\bigl|f(iy)\bigr|\leq \bigl|g(iy)\bigr|\quad \text{при  всех $y\in \RR$};
\end{equation}
\item\label{fgii} существуют ц.ф.э.т.
 $g\neq 0$ c последовательностью всех её нулей (с учётом кратности) $\Zero_g \supset {\sf W}$ и с ограничением ${\Zero_g}\cap \CC_{\rght} ={\sf W}$, а также ц.ф.э.т. $f\neq 0$ с $f({\sf Z})=0$, для которых выполнено \eqref{fgiR};
\item\label{fgiii} существует число  $C\in \RR$, с которым 
\begin{equation}\label{Zld}
l_{\sf Z}(r,R):=\sum_{r<{\sf z}_k\leq R} \frac{1}{{\sf z}_k} \leq l_{\sf W}(r,R) 
+C\quad \text{при всех\/ $0\leq r<R<+\infty$}.
\end{equation}
\end{enumerate}
\end{thMR}
В \cite{MR} установлены  и иные весьма содержательные вытекающие из теоремы Мальявена\,--\,Рубела  многочисленные результаты о росте ц.ф.э.т. вдоль прямой, о полноте экспоненциальных систем в пространствах функций, голоморфных в полосе, об  аналитическом продолжении и пр.

Характеристика $l_{\sf Z}(R)\overset{\eqref{Zld}}{:=}l_{\sf Z}(0,R)$ для 
положительных последовательностей ${\sf Z}$
называется в  \cite[(2.2)]{MR}, \cite[гл.~22, определение]{RC} {\it 
характеристическим логарифмом\/} последовательности $\sf Z$. 
Развивая её, в \cite[введение]{Kha88} частично, и детально в \cite{KhaD88}, \cite{Kha89}, \cite[\S~1]{kh91AA}, \cite[\S~0]{KhDD92}  (см. также \cite[3.2.1]{Khsur}) были 
 определены  \textit{ правый и левый характеристические логарифмы} последовательности ${\sf Z}$ в правой и левой открытых полуплоскостях 
$\CC_{\rght}$ и $\CC_{\lft}$  как
\begin{subequations}\label{logZC}
\begin{align}
l_{{\sf Z}}^{\rght}(R)&:=\sum_{\substack{0 < |{\sf z}_k|\leq R\\{\sf z}_k \in \CC_{\rght}}} \Re \frac{1}{{\sf z}_k}, \quad 0< R\leq +\infty, 
\tag{\ref{logZC}r}\label{df:dD+}\\
l_{{\sf Z}}^{\lft}(R)&:=
\sum_{\substack{0< |{\sf z}_k|\leq R\\{\sf z}_k \in \CC_{\lft}}} 
-\Re \frac{1}{{\sf z}_k}, \quad 0< R\leq +\infty,
\tag{\ref{logZC}l}\label{df:dD-}
\end{align}
\end{subequations}
а также   {\it правая\/} и {\it левая  логарифмические меры интервалов\/}  
$(r,R]\subset \RR^+$  
\begin{subequations}\label{df:l}
\begin{align}
l_{{\sf Z}}^{\rght}(r, R)&\overset{\eqref{df:dD+}}{:=}
l_{{\sf Z}}^{\rght}(R)-l_{{\sf Z}}^{\rght}(r),
\quad 0\leq r < R \leq +\infty,
\tag{\ref{df:l}r}\label{df:dDl+}\\
l_{{\sf Z}}^{\lft}(r, R)&\overset{\eqref{df:dD-}}{:=}l_{{\sf Z}}^{\lft}(R)-l_{{\sf Z}}^{\lft}(r), \quad 0\leq r < R \leq +\infty,
\tag{\ref{df:l}l}\label{df:dDl-}
\\
\intertext{которые порождают {\it логарифмическую субмеру  интервалов\/} для  $\sf Z$} 
l_{{\sf Z}}(r, R)&:=\max \bigl\{ l_{{\sf Z}}^{\lft}(r, R), l_{{\sf Z}}^{\rght}(r,
R)\bigr\}, \quad 0\leq r < R \leq +\infty .
\tag{\ref{df:l}m}\label{df:dDlL}
\end{align}
\end{subequations}
Для ${\sf Z}=\varnothing$ по определению $l_{\sf Z}(r,R)\equiv 0$ при всех 
$0\leq  r < R \leq +\infty$.

В  конце 1980-х гг.  в статьях  \cite[основная теорема ]{Kha88} (начальный  частный случай) и \cite[основная теорема]{KhaD88}, \cite[основная теорема]{Kha89} (см. также \cite[теорема 2.4.2]{KhDD92} и \cite[теорема 3.2.1]{Khsur}) были получены результаты без каких-либо  ограничений  на  последовательности   \eqref{Z}--\eqref{Zt}, кроме 
 конечности их  верхней плотности при порядке $1$ и ${\sf W}\subset \CC_{\overline \rght}$,  но при этом рассматривались ограничения более мягкие , чем  в \eqref{fgiR}, вида $\ln \bigl|f(iy)\bigr|\leq \ln \bigl|g(iy)\bigr|
+\varepsilon |y|$ при  всех $y\in \RR$, лежащих вне некоторого множества  $E\subset \RR$ конечной {\it лебеговой линейной меры\/} $\lambda_{\RR}$ на $\RR$  с произвольно малым $\varepsilon>0$. 

В начале 1990-х гг. в  \cite{kh91AA} удалось достичь уровня ограничения \eqref{fgiR} для последовательностей \textit{ комплексных\/} точек при некотором дополнительном условии на ${\sf Z}$. В \cite[(0.2)]{kh91AA} последовательность 
${\sf Z}\overset{\eqref{Z}}{\subset} \CC$ называется {\it отделённой\/}
(углами) {\it от мнимой оси,\/}  
 если для некоторого  числа $d>0$
\begin{equation}\label{con:dis}
|\Re {\sf z}_k|\geq d|{\sf z}_k|>0\quad\text{ при всех $k=1,2, \dots$}.
\end{equation}
Для функции $v\colon i\RR\to \RR_{\pm \infty}:=\{-\infty\}\cup \RR\cup \{+\infty\}$, локально интегрируемой по {\it линейной мере Лебега $\lambda_{i\RR}$ на\/}  $i\RR$,  определим {\it логарифмический интеграл\/}
\begin{equation}\label{fK:abp+}
J_{i\RR}(r,R; v):=
\frac{1}{2\pi}\int_r^R \frac{v(-iy)+v(iy)}{y^2} \dd y, \quad 0<r<R<+\infty,
\end{equation}
по всевозможным интервалам $(r,R] \subset \RR^+$ с левым концом $r>0$.

Теорема Kh  --- это  распространение  
 эквивалентности \eqref{fgii}$\Leftrightarrow$\eqref{fgiii} теоремы Мальявена\,--\,Рубела на  комплексные  последовательности \eqref{Z}--\eqref{Zt}.

\begin{theoKh}[{\rm \cite[основная теорема]{kh91AA}}]\label{theoKh1} 
Пусть комплексная  последоват\-е\-л\-ь\-ность ${\sf Z}\subset \CC$ и 
последовательность нулей $\Zero_g$ ц.ф.э.т. $g$ 
отделены углами от $i\RR$ в смысле \eqref{con:dis}. Тогда эквивалентны два утверждения 
\begin{enumerate}[{\rm (i)}]
\item\label{khi} существует ц.ф.э.т. $f\neq 0$ с $f({\sf Z})=0$, для которой выполнено
\eqref{fgiR};
 \item\label{khii} существует число  $C\in \RR$, для которого 
\begin{equation}\label{lJ}
l_{\sf Z}(r,R)\overset{\eqref{fK:abp+}}{\leq} J_{i\RR}\bigl(r,R; \ln |g|\bigr)+C
\quad\text{при всех $1\leq r<R<+\infty$.}
\end{equation}
\end{enumerate}
\end{theoKh}
Основные результаты  в  \cite{Kha88}, \cite{Kha89}, \cite{kh91AA}, 
доказываются ad ovo без использования  теоремы Мальявена\,--\,Рубела,  хотя, безусловно, ряд  важных  идей для их  доказательства   почерпнуты  из основополагающей работы П. Мальявена и Л.\,А. Рубела \cite{MR}. К таковым относится и ключевая идея выметания на мнимую ось, но уже не только меры или заряда с носителем на  $\RR^+$, но и мер и зарядов, <<размазанных>> по  всей  плоскости $\CC$.      

Важные и многочисленные примеры и контрпримеры по информативности и  пределах применимости логарифмических функций интервалов \eqref{df:l}, утверждения об их взаимосвязях с другими характеристиками последовательностей \eqref{Z}, а также точности результатов, связанных с ними, имеются  в \cite[теорема 8.1, 10]{MR}, \cite[\S~4, примеры]{Kha91e}, \cite[\S~7]{Kha91I}, \cite[\S~1, теоремы 3 и 6]{Kha09},
\cite[примеры 3.2.1, 3.2.2]{Khsur}.

\subsection{Развитие предшествующих результатов для ц.ф.э.т.}\label{ss:ef} Теоремы \ref{thMR+} и \ref{thKh2} этого подпараграфа будут выведены из соответствующих более общих результатов  для субгармонических функций в подпараграфе \ref{sver1-3} из \S~\ref{mainr}. 
 Следующая теорема \ref{thMR+} распространяет теорему Мальявена\,--\,Рубела  на  комплексные последовательности \eqref{Z}.   
\begin{theorem}\label{thMR+}
Пусть последовательность 
${\sf Z}$ из \eqref{Z}
отделена углами от мнимой  оси $i\RR$ в смысле \eqref{con:dis},  а ${\sf W}\subset \CC$ --- последовательность из \eqref{Z}--\eqref{Zt}, удовлетворяющая хотя бы одному из следующих трёх условий:
\begin{equation}\label{rlRZ}
\begin{split}
{\bf [r]}\quad  {\sf W}\subset \CC_{\overline \lft}
\qquad \text{или}  \qquad  {\bf [l]} \quad {\sf W}\subset \CC_{\overline \lft}\\
 \text{или} \quad  {\bf [R]}\quad  \sup_{r\geq 1}\left|\sum_{1<|{\sf w}_k|\leq r}\Re\frac{1}{{\sf w}_k}\right|<+\infty.
\end{split}
\end{equation} 

Тогда  эквивалентны  следующие три утверждения:
\begin{enumerate}[{\rm (i)}]
\item\label{fgiC} для любой ц.ф.э.т. $g\neq 0$, обращающейся в нуль на  $ {\sf W}$, и для  любого числа $P\in \RR^+$
 найдутся  ц.ф.э.т. $f\neq 0$, обращающаяся в нуль на  ${\sf Z}$, и 
борелевское подмножество $E\subset \RR$, для которых 
\begin{equation}\label{fEP}
\begin{split}
\bigl|f(iy)\bigr|&\leq \bigl|g(iy)\bigr|\quad \text{для каждого $y\in \RR\setminus E$, где}\\ 
\lambda_{\RR}&\bigl(E\setminus [-r,r] \bigr)
=O\Bigl(\frac{1}{r^P}\Bigr)\quad\text{при $r\to +\infty$.}
\end{split}
\end{equation}
\item\label{fgiiC} существуют ц.ф.э.т.
 $g\neq 0$, для которой  в случаях \eqref{rlRZ}  выполнено соответствующее требование  для $\Zero_g$ из следующих трёх:
\begin{equation}\label{rlRZ+}
\begin{split}
{\bf [r]}\quad \Zero_g\cap \CC_{\rght}={\sf W}\cap \CC_{\rght}
\qquad \text{или}  
\qquad  {\bf [l]}\quad \Zero_g\cap \CC_{\lft}={\sf W}\cap \CC_{\lft}\\
 \text{или} \quad  {\bf [R]}\quad  \text{либо  {\bf [r],} либо  {\bf [l]} из предшествующих}.
\end{split}
\end{equation}
и ц.ф.э.т. $f\neq 0$ с $f({\sf Z})=0$, для которых выполнено \eqref{fEP} при $P=0$;
\item\label{fgiiiC} существует число  $C\in \RR$, для которого  
\begin{equation}\label{ZldC}
l_{\sf Z}(r,R) \leq l_{\sf W}(r,R) 
+C\quad \text{при всех\/ $0\leq r<R<+\infty$}.
\end{equation}
\end{enumerate}
Если дополнительно предположить, что  последовательность ${\sf W}$ также отделена углами от $i\RR$, то  в \eqref{fEP} можно добиться  пустого $E:=\varnothing$, если утверждение \eqref{fgiC} 
заменить на утверждение 
\begin{enumerate}[{\rm (i$'$)}]
\item\label{fgiCS} для любой ц.ф.э.т. $g$ с последовательностью нулей $\Zero_g$,  отделённой углами от $i\RR$ и обращающейся в нуль на  ${\sf W}$,  найдётся ц.ф.э.т. $f\neq 0$, для которой
$f({\sf Z})=0$ и выполнено \eqref{fgiR}.
\end{enumerate}
\end{theorem}

Теорему Kh обобщает и развивает 
\begin{theorem}\label{thKh2} Пусть последовательность ${\sf Z}$ такая же, как в 
теореме  {\rm Kh}, но $g\neq 0$ --- произвольная ц.ф.э.т. Тогда 
утверждение \eqref{khii} теоремы  {\rm Kh} с соотношением \eqref{lJ} эквивалентно утверждению 
\begin{enumerate}[{\rm (i)}]
\item\label{fgiCSKh} для любого числа $P\in \RR^+$ найдётся  ц.ф.э.т. $f\neq 0$, обращающаяся в нуль на  ${\sf Z}$ и удовлетворяющая соотношениям  \eqref{fEP}.
\end{enumerate}
\end{theorem}

В настоящей статье не затрагиваются возможные применения наших основных результатов к вопросам нетривиальности весовых пространств ц.ф.э.т., полноты экспоненциальных систем в пространствах функций, к теоремам единственности для ц.ф.э.т., к существованию ц.ф.э.т.-мульт\-и\-п\-л\-и\-к\-а\-т\-о\-р\-ов, гасящих рост целой функции вдоль прямой, к представлению мероморфных функций в виде отношения ц.ф.э.т. с ограничениями на рост этих ц.ф.э.т. вдоль прямой, к аналитическому продолжению рядов, к задачам спектрального анализа-синтеза в пространствах голоморфных функций и пр. подобно   \cite{MR}--\cite{KhDD92}. Описания общей схемы и методики таких применений можно найти, например, в \cite{Khsur} и  \cite[гл.~2]{KhaRozKha19}.
Эти применения предполагается изложить в ином месте.   

\section{Основные результаты}\label{mainr}
\setcounter{equation}{0}

\subsection{Основные используемые определения и обозначения}\label{nod} 
\subsubsection{Меры и заряды}\label{1_2_4}
Как и в \cite{KhI}--\cite{KhII}, для борелевского подмножества $S\subset \CC$ через $\mathcal B (S)$  обозначаем класс всех \textit{борелевских подмножеств\/} в $S$, а через $\mathcal B_{\rm b} (S)\subset \mathcal B (S)$ --- класс всех {\it ограниченных в\/ $\CC$} борелевских подмножеств в $S$; 
$\Meas (S)$  --- класс всех {\it счетно-аддитивных функций на   $\mathcal B (S)$ со значениями на {\it расширенной числовой прямой\/} $\RR_{\pm\infty}$, конечных на  $\mathcal B_{\rm b} (S)\subset \mathcal B (S)$.\/}
Элементы из $\Meas (S)$ называем {\it зарядами\/} на $S\subset \CC$. Через  $\Meas^+ (S)\subset \Meas (S)$ обозначаем  подкласс всех {\it положительных зарядов\/} на $S$,  или просто \textit{мер\/} на $S$.  
Для заряда $\upnu$ через
$\upnu^+:=\sup\{0,\upnu\}$, $\upnu^-:=(-\upnu)^+$ и $|\upnu|:=\upnu^++\upnu^-$, обозначаем {\it верхнюю, нижнюю\/ {\rm и} полную вариации заряда\/
$\upnu=\upnu^+ -\upnu^-$.} Через $\updelta_z$ обозначаем {\it меру Дирака в точке $z\in \CC_{\infty}$,\/} т.е. вероятностную меру с носителем $\supp \updelta_z=\{z\}$. Заряд $\upnu \in \Meas(S)$ {\it целочисленный,} если для любого $B\in \mathcal B_{\rm b} (S)$ значение $\upnu(B)$ --- целое число.
Далее через  $D(z,r):=\{z' \in \CC \colon |z'-z|<r\}$,  
$\overline{D}(z,r):=\{z' \in \CC \colon |z'-z|\leq r\}$ и $\partial \overline{D}(z,r):=\overline{D}(z,r)\setminus{D}(z,r)$ обозначаем соответственно
{\it открытый круг},   {\it замкнутый круг} и {\it окружность  с центром\/ $z\in \CC$ радиуса\/ $r\in \RR^+$}. В частности,  
 $D(r):=D(0,r)$,  $\overline D(r):=\overline D(0,r)$, а   $\DD:=D(1)$ и $\overline \DD:=\overline D(1)$ --- соответственно {\it открытый\/} и {\it замкнутый 
единичный круг.\/} Для заряда $\upnu\in \Meas (S )$ и круга $\overline D (z,r)\subset S$ полагаем
\begin{equation}\label{df:nup} 
\upnu (z,r):=\upnu \bigl(\,\overline D(z,r)\bigr),\quad \upnu^{\rad}(r):=\upnu(0,r)=\upnu\bigl(\overline D(r)\bigr) 
\end{equation} 
--- {\it считающие  функция с
центром\/ $z$} и {\it радиальная считающая функция} заряда $\upnu$. 
Последовательности ${\sf Z}$ из \eqref{Z} сопоставляем 
 {\it считающую меру}
\begin{equation}\label{df:divmn}
n_{\sf Z}:=\sum_k \updelta_{{\sf z}_k}, \quad  
n_{\sf Z} (S)= \sum_{{\sf z}_k\in S} 1 \in \{0\}\cup \NN\cup \{+\infty\}
 \quad \text{для $S\subset \CC$}.
\end{equation}
Определения   \eqref{df:nup}--\eqref{df:divmn}  согласованы  с  \eqref{nZrad},
поскольку $(n_{\sf Z})^{\rad}=n_{\sf Z}^{\rad}$. 

Две последовательности  из \eqref{Z} равны или совпадают,  если  $n_{\sf Z}\overset{\eqref{df:divmn}}{=} n_{\sf W}$ как меры. 
При этом носитель $\supp {\sf Z}$ --- это $\supp n_{\sf Z}$; $z \in {\sf Z}$ или $z \notin {\sf Z}$ соответственно означает, что $z \in \supp {\sf Z}$ или  $z \notin \supp {\sf Z}$. Для подмножества $S \subset {\CC}$ запись ${\sf Z} \subset S$ означает, что 
$\supp {\sf Z} \subset S$; ${\sf Z} \cap S$ --- сужение последовательности  $\sf Z$ на $S$ с  считающей мерой $n_{\sf Z}\bigm|_S$, суженной на $S$.
Последовательность точек ${\sf Z}$ включена (содержится) в ${\sf W}$, если  $n_{\sf Z}\leq n_{\sf W}$. При этом пишем ${\sf W}\subset {\sf Z}$ и говорим, что ${\sf W}$ --- {\it подпоследовательность\/} из ${\sf Z}$.
 
{\it Объединение\/}  ${\sf Z} \cup {\sf W}$
через считающие меры  задается равенством 
$n_{{\sf Z} \cup {\sf W}}=n_{\sf Z}+n_{\sf W}$, а \textit{пересечение\/} ${\sf Z} \cap {\sf W}$ через считающую меру 
$n_{{\sf Z} \cap {\sf W}}=\inf \{n_{\sf Z},n_{\sf W}\}$. 

На последовательностях  точек операции и отношения, отличные от приведенных выше,  понимаются поэлементно. Так,
$z{\sf Z} :=\{ z{\sf z}_k\}$ и т.п.

\subsubsection{Логарифмические характеристики  мер и зарядов} 
Распространим определения правой и левой логарифмических функций интервалов для ${\sf Z}\subset \CC$  из  \eqref{df:l}  на произвольные заряды  $\upnu\in \Meas (\CC)$:
\begin{subequations}\label{df:lm}
\begin{align}
l_{\upnu}^{\rght}(r, R)&\overset{\eqref{df:dDl+}}{:=}\int_{\substack{r < | z|\leq R\\ \Re z >0}} \Re \frac{1}{ z} \dd \upnu(z), \quad 0< r < R \leq +\infty ,
\tag{\ref{df:lm}r}\label{df:dDlm+}\\
l_{\upnu}^{\lft}(r, R)&\overset{\eqref{df:dDl-}}{:=}\int_{\substack{r< |z|\leq R\\ 
\Re z<0}}\Re \Bigl(-\frac{1}{ z}\Bigr) \dd \upnu(z) ,  \quad 0< r < R \leq  +\infty ,
\tag{\ref{df:lm}l}\label{df:dDlm-}\\
\intertext{--- соответственно {\it правая\/} и {\it левая логарифмические функции интервалов\/} для заряда $\upnu$,  которые в случае  {\it меры\/} $\upmu\in \Meas^+(\CC)$ порождают
{\it правую\/} $l_{\upmu}^{\rght}$ и {\it левую $l_{\upmu}^{\lft}$ логарифмические меры интервалов\/} для меры $\upmu$, а также 
 {\it логарифмическую субмеру интервалов\/} $l_{\upmu}$ для $\upmu$, а именно:}
l_{\upmu}(r, R)&\overset{\eqref{df:dDlL}}{:=}\max \bigl\{ l_{\upmu}^{\lft}(r, R), l_{\upmu}^{\rght}(r,R)\bigr\}, \quad 0< r < R \leq +\infty. 
\tag{\ref{df:lm}m}\label{df:dDlLm}
\end{align}
\end{subequations}
Если  $0\notin \supp \upnu$ или, более общ\'о, заряд $\upnu$ принадлежит классу сходимости при порядке $1$ около нуля \cite[2.2, предложение 2.2]{KhI}, т.\,е.
\begin{equation*}
\int_0^{r_0}\frac{|\upnu|^{\rad}(t)}{t^2}\dd t<+\infty\quad\text{при некотором $0<r_0\in \RR^+$},
\end{equation*}
то  по непрерывности при $0<r\to 0$ определены и введённые  ранее в \cite[(2.2)]{MR}, \cite[гл.~22, Определение]{RC} лишь для положительных последовательностей ${\sf Z}$ \textit{правый\/} и {\it левый характеристические логарифмы заряда\/} $\upnu$:
\begin{subequations}\label{df:dDlL1}
\begin{align}
l_{\upnu}^{\rght}(R)&\overset{\eqref{df:dD+}}{:=}l_{\upnu}^{\rght}(0,R):=
\lim_{0<r\to 0}l_{\upnu}^{\rght}(r,R),
\tag{\ref{df:dDlL1}r}\label{{df:dDlL1}r}
\\ 
l_{\upnu}^{\lft}(R)&\overset{\eqref{df:dD-}}{:=}l_{\upnu}^{\lft}(0,R)
:=\lim_{0<r\to 0}l_{\upnu}^{\lft}(r,R), 
\tag{\ref{df:dDlL1}l}\label{{df:dDlL1}l}\\
\intertext{а в случае {\it положительной меры\/} $\upmu\in \Meas^+(\CC)$ --- двусторонний {\it характеристический логарифм меры\/} $\upmu$:}
l_{\upmu}(R)&:=l_{\upmu}(0,R)
:=\lim_{0<r\to 0}l_{\upmu}(r,R).
\tag{\ref{df:dDlL1}m}\label{{df:dDlL1}m}
\end{align}
\end{subequations}
Для пары мер $\upnu, \upmu \in \Meas^+(\CC)$ ключевую роль на протяжении всей статьи будет играть следующее  соотношение:
\begin{equation}\label{lJMmul}
\boxed{\sup_{1\leq r<R<+\infty}
\Bigl(l_{\upnu}(r,R)-l_{\upmu}(r,R)\Bigr)<+\infty.}
\end{equation} 
\begin{definition}\label{ab} Мера\/ $\upmu \in \Meas^+(\CC)$  {\it  логарифмически доминирует над мерой\/ $\upnu \in \Meas^+(\CC)$,} или {\it мажорирует\/} $\upnu \in \Meas^+(\CC)$,  {\it вдоль вещественной оси $\RR$,} или {\it в направлениях $0$ и $\pi$,} если выполнено соотношение 
\eqref{lJMmul}, и при этом пишем $\upnu\curlyeqprec_{|\Re|} \upmu$. В частности, для двух последовательностей ${\sf Z}$ и ${\sf W}$ из \eqref{Z} пишем ${\sf Z}\curlyeqprec_{|\Re|}{\sf W}$, если $n_{\sf Z}\curlyeqprec_{|\Re|}n_{\sf W}$ для считающих мер \eqref{df:divmn}. 
\end{definition}
 \begin{remark}
Это определение для предпорядка $\curlyeqprec_{|\Re|}$ на  мерах или последовательностях точек полностью согласовано с введённым и исследованным в \cite[3, лемма 3.1]{MR} отношением    ${\sf Z}\prec {\sf W}$ и в \cite[гл.~22]{RC} ${\sf Z} <{\sf W}$ лишь  на   \textit{положительных\/}  последовательностях
${\sf Z},{\sf W}\subset \RR^+$. Отношение предпорядка  $\curlyeqprec_{|\Re|}$
в подобном обозначении было использовано в \cite[определение 7.1]{KhaRozKha19} для понятия {\it аффинного выметания.\/} В терминологии  
\cite[определения 5.2, 7.1]{KhaRozKha19},
\cite[определение 1]{KhaKha19}, \cite[определение 1]{MenKha20}  приведённое здесь отношение $\upnu\curlyeqprec_{|\Re|} \upmu$ означает, что мера $\upmu$ --- это аффинное выметание  меры $\upnu$ относительно  класса функций $z\underset{\text{\rm \tiny $0\neq z\in \CC$}}{\longmapsto} \bigl|\Re (1/z)\bigr|$, продолженных нулём   вне  всевозможных концентрических колец $\bigl\{z\in \CC\colon r<|z|\leq R \bigr\}$.
\end{remark}

Для величин $-\infty<\alpha< \beta \leq \alpha +2\pi<+\infty$ определяем 
{\it замкнутый угол с вершиной в нуле $\angle [\alpha, \beta]:=\{z\in \CC\colon \alpha\leq \arg z\leq \beta\}$ с вершиной в нуле.}

\textit{Заряд $\upnu \in \Meas(\CC)$ отделён углами от мнимой оси\/} $i\RR$, если 
для некоторого числа $\varepsilon>0$ его носитель $\supp \upnu$
не пересекается с объединением 
\begin{equation}\label{anglei}
\angle[\pi/2-\varepsilon, \pi/2+\varepsilon]\bigcup \angle[-\pi/2-\varepsilon, -\pi/2+\varepsilon], \quad \varepsilon >0,   
\end{equation} 
пары замкнутых вертикальных углов раствора $2\varepsilon>0$ с вершиной в нуле и с биссектрисой $i\RR$. 
В частности, \textit{комплексная  последовательность\/ ${\sf Z}\subset \CC$ отделена  углами от мнимой оси\/} $i\RR$, если такова её считающей мера $n_{\sf Z}$, что эквивалентно требованию  \eqref{con:dis}.

\subsubsection{Функции}
Пусть $S\subset \CC$. Классы $\Hol(S)$ и  $\sbh (S)$ состоят из сужений на $S$ функций, соответственно  {\it голоморфных\/} и {\it субгармонических\/}  в каком-либо открытом множестве,  включающем в себя $S$ \cite{Rans}, \cite{HK};  
$\har (S):=\sbh (S)\cap \bigl(-\sbh(S)\bigr)$ --- класс {\it гармонических\/} функций на $S$; $\dsbh (S):=\sbh (S) -\sbh(S)$ --- класс \textit{$\updelta$-субгармонических,\/} или разностей субгармонических,  \textit{функций\/} на $S$ \cite{Ar_d}, \cite{Arsove53p},  \cite[2]{Gr}, \cite[3.1]{KhaRoz18}.
Полагаем  $\Hol_*(S):=\Hol(S)\setminus \{0\}$.  Тождественную $-\infty$ или $+\infty$ на $S$ обозначаем соответственно $\boldsymbol{-\infty}\in \sbh (S)$
или $\boldsymbol{+\infty}\in -\sbh (S)\subset \dsbh (S)$; $\sbh_*(S):=\sbh (S)\setminus \{\boldsymbol{-\infty}\}$, $\dsbh_*(S):=\dsbh (S)\setminus \{\boldsymbol{\pm\infty}\}$. 

Пусть $r_0\in \RR$. {\it Тип\/}  функции $f\colon [r_0,+\infty)\to \RR_{\pm \infty}$ {\it при порядке\/} $p\in \RR^+$ (около $\infty$) в обозначении $f^+:=\sup\{f,0\}$ определяется как  \cite[2.1]{KhI}
\begin{equation}\label{typevf}
\type_p[f]\overset{\text{\cite[(2.1t)]{KhI}}}{:=}
\type_p^{\infty}[f]:=
\limsup_{r\to +\infty}
\frac{f^+(r)}{r^p}\in \RR^+\cup \{+\infty\}.
\end{equation}
{\it Порядок\/} этой  функции $f$ (около $\infty$) определяется как 
\cite[2.1]{KhI} 
\begin{equation}\label{ord}
\ord[f]\overset{\text{\cite[(2.1a)]{KhI}}}{:=}\ord_{\infty}[f]
:=\limsup\frac{\ln \bigl(1+f^+(r)\bigr)}{\ln r}\in \RR^+\cup \{+\infty\}. 
\end{equation} 
Для $z\in \CC$, $r \in \RR^+$ определим  интегральное средние по   окружности $\partial \overline D(z, r)$  
от интегрируемой на ней   фу\-н\-к\-ц\-ии $ v\colon \partial 
  D(z,r)\to \RR_{\pm\infty}$:
\begin{subequations}\label{df:MCB}
\begin{align}
\mathsf{C}_v(z, r)&
\:=\frac{1}{2\pi} \int_{0}^{2\pi}  v(z+re^{i\theta}) \dd \theta, \quad \mathsf{C}_v(r):=\mathsf{C}_v(0, r),\tag{\ref{df:MCB}C}\label{df:MCBc}
\\
\intertext{а также верхнюю грань функции $v\colon \partial D(z,r)\to \RR_{\pm \infty}$ на $\partial D(z,r)$:}
 \mathsf{M}_v(z,r)&:=\sup_{z'\in \partial D(z,r)}v(z'), 
\quad \mathsf{M}_v(r):=\mathsf{M}_v(0, r),  
\tag{\ref{df:MCB}M}\label{df:MCBm}
\end{align}
\end{subequations}
которая для функции  $v\in \sbh \bigl( \overline D(z, r)\bigr)$ совпадает с верхней гранью 
функции $v$ в круге $\overline D(z, r)$  \cite[определение 2.6.7]{Rans}, \cite{HK}. Для  функции $v\in \sbh(\CC)$ 
её {\it порядок\/} и {\it тип при порядке\/} $p$ определяются как  
\cite[замечание 2.1]{KhI}
\begin{equation}\label{typev}
\ord[v] \overset{\eqref{df:MCBm}}{:=}\ord[{\mathsf{M}_{v}}], \quad 
\type_p[v] \overset{\eqref{df:MCBm}}{:=}\type_p[{\mathsf{M}_{v}}].
\end{equation}

Заряд $\upnu\in \Meas (\CC)$ \textit{конечной верхней плотности\/} при порядке $1$, или \textit{конечного типа\/} при порядке $1$, если \cite[\S~4]{KhI}
\begin{equation}\label{fden}
{\type_1}[\upnu]\overset{\eqref{typevf}}{:=}\type_1\bigl [|\upnu|^{\rad}\bigr]:= \limsup_{r\to+\infty}\frac{|\upnu|^{\rad}(r)}{r}<+\infty.
\end{equation}

Функция $v\in \dsbh_*(\CC) $ конечного тип, если она представима
в виде  разности двух функций из $\sbh_*(\CC)$ конечного типа при порядке $1$.

Функция $f\in \Hol(\CC)$ ---
{\it целая функция\/  экспоненциального типа\/} (пишем {\it ц.ф.э.т.\/}), если 
 $\ln |f|\in \sbh(\CC)$ конечного типа при порядке $1$.

Для открытого множества $O\subset \CC$ меру Рисса функции ${u}\in \sbh_* ( O)$ обозначаем как 
$\frac{1}{2\pi}\Delta {u}\in \mathcal M ^+(O)$, где оператор Лапласа $\Delta$ действует в смысле теории обобщённых функций. 
Для функции  $\boldsymbol{-\infty}$ её мера Рисса по
определению равна $+\infty$ на любом подмножестве из $ O$.

Взаимосвязь между распределением нулей целых  функций и мерой Ри\-с\-са  субгармонических функций определяет следующее 
\begin{propos}[{\rm \cite[лемма 10.10]{Hayman}, \cite[теорема 3.7.8]{Rans}}]\label{prfZ} Если  $f\neq 0$ --- целая функция, то считающая мера последовательности её нулей      
\begin{equation}\label{ZS}
n_{\Zero_f}\overset{\eqref{df:divmn}}{=}\frac{1}{2\pi}\Delta_{\ln|f|}\in \Meas^+(\CC). 
\end{equation}
--- целочисленная мера. Обратно, если мера $\upnu \in \Meas^+(\CC)$ целочисленная, то для любой функции $u\in \sbh_*(\CC)$ с мерой Рисса 
$\frac{1}{2\pi}\Delta_u=\upnu$ существует целая функция  $f\in \Hol_*(\CC)$ со считающей функцией нулей $n_{\Zero_f}=\upnu$, для которой $\ln |f|=u$ на всей плоскости $\CC$.  
\end{propos}

Функция $f\in \Hol_*(\CC)$ с  {\it последовательностью\/} всех её  {\it нулей,\/}
или корней, $\Zero_f$,  перенумерованной  с учетом кратности {\it обращается в нуль на\/} ${\sf Z}$, если  ${\sf Z} \subset \Zero_f$. При этом пишем $f({\sf Z})=0$.

\subsection{Необходимые условия для мер Рисса}\label{mT}
Следующий результат в гораздо более общей форме установлен  
 в \cite{Kha20}. 

\begin{propos}[{\cite[теорема 1]{Kha20}}]\label{th12} Пусть функция $M\in \sbh_*(\CC)$ конечного типа $\type_1[M]<+\infty$ при порядке $1$ с мерой Рисса $\upmu:=\frac{1}{2\pi}\Delta M$ и 
для функции $u\in \sbh_*(\CC)$ конечного типа $\type_1[u]<+\infty$ при порядке $1$ с мерой Рисса $\frac{1}{2\pi}\Delta u\geq \upnu\in \Meas^+(\CC)$ имеет место ограничение 
\begin{equation*}
u(iy)+u(-iy)\leq M(iy)+M(-iy)
\quad \text{при всех $y\in \RR^+\setminus E$},
\end{equation*}
где $E\in \mathcal B( \RR^+)$
конечной меры $\lambda_{\RR}(E)<+\infty$. Тогда 
\begin{equation*}
\sup_{1\leq r<R<+\infty}
\biggl(l_{\upnu}(r,R)-\min\Bigl\{l_{\upmu}^{\rght}(r,R),
l_{\upmu}^{\lft}(r,R), J_{i\RR}(r,R; M)\Bigl\}\biggr)<+\infty
\end{equation*}
и тем более, по определению \eqref{df:dDlLm},
имеем \eqref{lJMmul}, т.е. $\upnu\curlyeqprec_{|\Re|} \upmu$.
\end{propos}

Взаимосвязь логарифмических характеристик мер и последовательностей  
с  логарифмическими  интегралами  \eqref{fK:abp+}
 по  интервалам $(r,R] \subset \RR^+$ устанавливает следующее   
\begin{propos}[{\cite[предложение 4.1, (4.19)]{KhII}}]\label{lemJl} 
Для  любой субгармонической функции $M\in \sbh_*(\CC)$ конечного типа $\type_1[M]<+\infty$ при порядке $1$ с мерой Рисса $\upmu :=\frac{1}{2\pi}\Delta M\in \Meas^+(\CC)$ имеют место соотношения
\begin{subequations}\label{Jll}
\begin{align}
\sup_{1\leq r<R<+\infty} &\Bigl|J_{i\RR}(r,R;M)-l_{\upmu}^{\rght}(r,R)\Bigr|
<+\infty,
\tag{\ref{Jll}r}\label{{Jll}r}
 \\
\sup_{1\leq r<R<+\infty} &\Bigl|J_{i\RR}(r,R;M)-l_{\upmu}^{\lft}(r,R)\Bigr|
<+\infty,
\tag{\ref{Jll}l}\label{{Jll}l}
\\
\sup_{1\leq r<R<+\infty} &\Bigl|J_{i\RR}(r,R;M)-l_{\upmu}(r,R)\Bigr|
<+\infty.
\tag{\ref{Jll}m}\label{{Jll}m} 
\end{align}
\end{subequations}
\end{propos}

\subsection{Формулировка основной теоремы}\label{MT} 
Далее неоднократно и зачастую без ссылки будут использоваться две классические теоремы о представлении целых и субгармонических функций конечного порядка. Мы формулируем их лишь в той минимальной форме, которая необходима для использования в данной статье. 
Последовательность ${\sf Z}$ из \eqref{Z} удовлетворяет {\it условию Линделёфа\/}   рода $1$, если 
\begin{equation}\label{con:LpZ} 
\sup_{r\geq 1}\Biggl|\sum_{1< |{\sf z}_k|\leq r} \frac{1}{{\sf z}_k}
\Biggr|<+\infty. 
\end{equation}

\begin{thWHL}[{\cite{Levin56}, \cite{Levin96}, \cite[2.10]{Boas}}]\label{pr:repef}
Если  $f\in \Hol_*(\CC)$ --- ц.ф.э.т., то последовательность её нулей  $\Zero_f$ конечной верхней плотности  $\type_1[\Zero_f]\in \RR^+$ и удовлетворяет условию  Линделёфа  \eqref{con:LpZ}.
Обратно, если  ${\sf Z}\subset \CC$ --- последовательность конечной верхней плотности  $\type_1[{\sf Z}]\in \RR^+$,  то  для любой целой функции $f\in \Hol_*(\CC)$ с  $\Zero_f={\sf Z}$ имеем представление 
\begin{equation}\label{repruZ}
\begin{split}
f=f_{\sf Z}e^h, \quad f_{\sf Z}&\in\Hol_*(\CC), 
\quad\Zero_{f_{\sf Z}}={\sf Z}, \quad h\in \Hol(\CC), \\
\limsup_{z\to \infty}&\frac{\ln |f_{\sf Z}|}{|z|\ln |z|}<+\infty ,
\quad \ord\bigl[\ln |f|\bigr]\leq 1.
\end{split}
\end{equation}
 При дополнительном условии Линделёфа \eqref{con:LpZ} рода $1$ для ${\sf Z}$  в представлении \eqref{repruZ}  можно  выбрать уже  ц.ф.э.т. $f_{\sf Z}$, а если $f$ нулевого типа $\type_2[f]=0$ при  порядке $2$, то можно даже  положить $h=0$.  
\end{thWHL}

Заряд $\upmu\in \Meas(\CC)$ удовлетворяет  {\it условию Линделёфа\/} рода $1$  \cite[\S~2]{kh91AA}, \cite[определение 2.3.1]{KhDD92}, \cite[определение 4.7]{KhI}, если  
\begin{equation}\label{con:Lp+} 
\sup_{r\geq 1}\Biggl|\int_{D(r)\setminus \overline \DD} \frac{1}{z}
\dd \upmu(z) \Biggr|<+\infty. 
\end{equation}

\begin{thWHLB}[{\cite[3, Теорема  12]{Arsove53p},  \cite[4.1, 4.2]{HK}, \cite[2.9.3]{Az}, \cite[6.1]{KhI}}]\label{pr:rep}
Если   $u\in \sbh_*(\CC)$ --- функция с мерой Рисса  $\upmu:=\frac{1}{2\pi}\Delta_u$
и  $\type_1[u]<+\infty$, то $\type_1[\upmu]<+\infty$  и мера $\upmu$ удовлетворяет условию Линделёфа \eqref{con:Lp+}. 
Обратно, если мера $\upmu\in \Meas^+(\CC)$ конечного типа  $\type_1[\upmu]<+\infty$ при порядке $1$, а   $u\in \sbh_*(\CC)$ --- функция с мерой Рисса  $\frac{1}{2\pi}\Delta_u=\upmu$, то имеет место представление
\begin{equation}\label{repru}
\begin{split}
u=&u_{\upmu}+H, \quad  u_{\upmu}\in \sbh_*(\CC),
 \quad   H\in \har (\CC), \quad \frac{1}{2\pi}\varDelta_{u_\upmu}=\upmu, \\ 
&\limsup_{z\to \infty}\frac{u_{\upmu}(z)}{|z|\ln |z|}<+\infty, \quad 
 \type_1[{\sf C}_u]<+\infty,
\end{split}
\end{equation}
где при условии\/ $\type_2[u]=0$ функция   $H$ ---  гармонический многочлен степени  $\deg H\leq 1$,  а если\/ $\upmu$ удовлетворяет условию Линделёфа \eqref{con:Lp+}  рода $1$, то в качестве функции $u_{\upmu}$ можно выбрать функцию с $\type_1[u_{\upmu}]<+\infty$.
\end{thWHLB}

\begin{maintheorem}\label{th21} Пусть меры $\upnu, \upmu\in \Meas^+(\CC)$ конечной верхней плотности  $\type_1[\upnu+\upmu]<+\infty$, мера $\upnu$ отделена  углами  от мнимой оси и выполнено условие  \eqref{lJMmul}, т.е. $\upnu\curlyeqprec_{|\Re|} \upmu$. 
Тогда для любой субгармонической функции   $M\in \sbh_*(\CC)$  конечного  типа  $\type_1[M]<+\infty$  с мерой Рисса $\frac{1}{2\pi}\Delta M\geq \upmu$ 
выполнены два утверждения: 
\begin{enumerate}
\item[{\rm [sbh]}] Найдётся  функция $U\in \sbh_*(\CC)$ конечного типа $\type_1[U]<+\infty$ с мерой Рисса $\frac{1}{2\pi}\Delta U\geq \upnu$,  для которой   
\begin{equation}\label{u=M}
U(iy)= M(iy)\quad\text{для всех $y\in \RR$.} 
\end{equation} 
\item[{\rm [Hol]}]
Для любого числа $P\in \RR^+$ и для любой  функции $u\in \sbh_*(\CC)$ с мерой Рисса $\frac{1}{2\pi}\Delta u=\upnu$ найдётся такая целая функция $h\in \Hol_*(\CC)$, что субгармоническая сумма  $u+\ln |h|\in \sbh_*(\CC)$ --- функция конечного типа $\type_1\bigl[u+\ln |h|\bigr]<+\infty$ при порядке $1$ и 
\begin{subequations}\label{ME}
\begin{align}
u(iy)+\ln \bigl|h(iy)\bigr|\leq M(iy) \quad\text{при всех $y\in \RR\setminus E$,}
\tag{\ref{ME}M}\label{vhM}
 \intertext{где 
$E\in \mathcal B (\RR)$ удовлетворяет условию} 
\lambda_{\RR}\bigl(E\setminus [-r,r]\bigr)=O\Bigl(\frac{1}{r^P}\Bigr)
\quad \text{при $r\to +\infty$}. 
\tag{\ref{ME}E}\label{PE}
\end{align}
\end{subequations}
\end{enumerate}
Если в дополнение при некотором $\varepsilon >0$ функция $M$ гармоническая на объединении  \eqref{anglei}  пары замкнутых вертикальных углов, содержащих  $i\RR$, то в заключении\/ {\rm [Hol]}   в \eqref{ME}  можно выбрать пустое $E:=\varnothing$. 
\end{maintheorem}
Доказательство основной теоремы   завершается  в последнем \S~\ref{proof:2}.

\subsection{Субгармонические версии теорем \ref{thMR+} и \ref{thKh2}}\label{sver1-3} 

Субгармоническая версия теоремы Мальявена\,--\,Рубела и теоремы \ref{thMR+}   --- следующая 

\begin{theorem}\label{thMR+m}
Пусть мера $\upnu\in \Meas^+(\CC)$ отделена углами от мнимой оси  $i\RR$, 
а  мера $ \upmu \in \Meas^+(\CC)$ конечной верхней плотности  
при порядке $1$, удовлетворяющая хотя бы одному из следующих трёх условий:
 \begin{equation}\label{rlR}
\begin{split}
{\bf [r]}\quad  \supp \upmu\subset \CC_{\overline \lft}
\qquad \text{или}  \qquad  {\bf [l]} \quad \supp \upmu\subset \CC_{\overline \lft}\\
 \text{или} \quad  {\bf [R]}\quad  \sup_{r\geq 1}\left|\int_{1<|z|\leq r}\Re\frac{1}{z}\dd \upmu(z)\right|<+\infty.
\end{split}
\end{equation}
Тогда  эквивалентны  следующие пять  утверждений:
\begin{enumerate}[{\rm (i)}]
\item[{\rm (0)}]\label{fgiC0} для любой функции $M\in \sbh_*(\CC)$ с $\type_1[M]<+\infty$ и с мерой Рисса $\frac{1}{2\pi}\Delta M\geq \upmu$ найдётся такая функция $U\in \sbh_*(\CC)$ с $\type_1[U]<+\infty$ и с мерой Рисса $\frac{1}{2\pi}\Delta U\geq \upnu$, что выполнено \eqref{u=M};

\item\label{fgiCs} для любых  функций $M\in \sbh_*(\CC)$ с $\type_1[M]<+\infty$  и с  мерой Рисса $\frac{1}{2\pi}\Delta M\geq \upmu$, а также   $u\in \sbh_*(\CC)$  с мерой Рисса 
$\frac{1}{2\pi}\Delta u= \upnu$ при любом $P\in \RR^+$ 
найдутся целая функция $h\in \Hol_*(\CC)$
и $E\in \mathcal B( \RR)$, для которых $\type_1\bigl[u+\ln|h|\bigr]<+\infty$ и выполнено \eqref{ME}; 

\item\label{fgiiCs} существуют функции  $M\in \sbh_*(\CC)$ с $\type_1[M]<+\infty$ и мерой Рисса $\upmu':=\frac{1}{2\pi}\Delta M\geq \upmu $, для которой выполнено соответствующее условие   
\begin{equation}\label{rlRZ+mu}
\begin{split}
&{\bf [r]}\quad \CC_{\rght}\cap \supp \upmu'
=\CC_{\rght}\cap \supp \upmu\\
\text{или}  
\quad  &{\bf [l]}\quad \CC_{\lft}\cap \supp \upmu'
=\CC_{\lft}\cap \supp \upmu\\
 \text{или} \quad  &{\bf [R]}\quad  \text{либо  {\bf [r],} либо  {\bf [l]} из предшествующих},
\end{split}
\end{equation}
а также $U\in \sbh_*(\CC)$  с $\type_1[U]<+\infty$  и мерой Рисса $\frac{1}{2\pi}\Delta U\geq \upnu$, для которых  $U(iy)\leq M(iy)$ для каждого  $y\in \RR\setminus E$, где $\lambda_{\RR}(E)<+\infty$;
\item\label{fgiiiCs} $\upnu\curlyeqprec_{|\Re|}\upmu$, т.е. выполнено \eqref{lJMmul};
\item\label{fgivCs} существует  строго возрастающая  последовательность  чисел $(r_n)_{n\in \NN}\subset \RR^+$, для которой 
$\lim\limits_{n\to\infty}{r_{n+1}}/{r_n}<+\infty$ и   
\begin{equation}\label{cprec}
 \limsup_{N\to  \infty}\sup\limits_{n\leq N}
\Bigl(l_{\upnu}(r_n,r_N)-l_{\upmu}(r_n,r_N)\Bigr)<+\infty.
\end{equation}
\end{enumerate}
Если дополнительно предположить, что  мера $\upmu$ также отделена углами от $i\RR$, то  в утверждении \eqref{fgiCs} для   \eqref{ME} можно добиться  пустого подмножества $E\subset \RR$, если утверждение \eqref{fgiCs} 
переписать в виде  
\begin{enumerate}[{\rm (i$'$)}]
\item\label{fgiCS+} 
 для любой функции $M\in \sbh_*(\CC)$ с $\type_1[M]<+\infty$ с мерой Рисса $\frac{1}{2\pi}\Delta M\geq \upmu$,  отделённой углами от $i\RR$,  и любой функции $u\in \sbh_*(\CC)$  с мерой Рисса 
$\frac{1}{2\pi}\Delta u= \upnu$ найдётся целая функция $h\in \Hol_*(\CC)$, для которой  $\type_1\bigl[u+\ln|h|\bigr]<+\infty$ и  \eqref{vhM}  выполнено для $E:=\varnothing$.
\end{enumerate}
\end{theorem}
\begin{proof} Доказательство проводится  по схеме 
\begin{equation}\label{diag}
\begin{array}{ccccc}
(0) &\longrightarrow & \text{\eqref{fgiiCs}} &\longleftarrow&\text{\eqref{fgiCs}}\\
 &  \nwarrow&  \downarrow&\nearrow&\\
\text{\eqref{fgivCs}} &  \longleftrightarrow& \text{\eqref{fgiiiCs}}&&\\
\end{array}
\end{equation}
с отдельным обсуждением дополнения после  \eqref{fgivCs}. 

Далее пишем  A$\overset{{\text{\it \tiny proof}}}{\Longrightarrow}$B
или A$\overset{\text{\it \tiny proof}}{\Longleftrightarrow}$B, если после этого доказывается или обсуждается  соответственно 
импликация  A$\Rightarrow$B  или эквивалентность A$\Leftrightarrow$B. 

(0)$\overset{{\text{\it \tiny proof}}}{\Longrightarrow}$\eqref{fgiiCs}$\overset{{\text{\it \tiny proof}}}{\Longleftarrow}$\eqref{fgiCs}. Для доказательства обеих этих  импликаций достаточно построить  функцию  $M\in \sbh_*(\CC)$ конечного типа  $\type_1[M]<+\infty$ с мерой Рисса $\upmu':=\frac{1}{2\pi}\Delta M\geq \upmu $, удовлетворяющей соотношению  $\upmu'\curlyeqprec_{|\Re|}\upmu$. Согласно теореме  Вейерштрасса\,--\,Адамара\,--\,Линделёфа\,--\,Брело такая функция найдётся, если для меры $\upmu$ можно построить меру $\upmu'\geq \upmu$, удовлетворяющую условию Линделёфа  \eqref{con:Lp+}  рода $1$  и одновременно  соотношению  $\upmu'\curlyeqprec_{|\Re|}\upmu$. 
Этот факт установлен ниже в предложении \ref{prL} с \eqref{mus}.

\eqref{fgiiCs}$\overset{{\text{\it \tiny proof}}}{\Longrightarrow}$\eqref{fgiiiCs}.
Эта импликация сразу следует из предложения \ref{th12}.

\eqref{fgiiiCs}$\overset{{\text{\it \tiny proof}}}{\Longrightarrow}$(0).
Это  частный случай части [sbh] основной теоремы.

\eqref{fgiiiCs}$\overset{{\text{\it \tiny proof}}}{\Longrightarrow}$\eqref{fgiCs}. Это частный случай части [Hol] основной теоремы.

\eqref{fgivCs}$\overset{{\text{\it \tiny proof}}}{\Longleftrightarrow}$\eqref{fgiiiCs}. Это предложение \ref{pr:rn}, основанное на свойствах логарифмических характеристик зарядов, исследованных в подпараграфе \ref{2_3_1}.  

Дополнение с утверждением (\ref{fgiCS+}$'$)  для отделённой углами от $i\RR$  меры $\upmu$ 
сразу следует из завершающего основную теорему утверждения. 
\end{proof}
\begin{proof}[Доказательство теоремы \ref{thMR+}]
Положим $\upnu:=n_{\sf Z}$, $\upmu:=n_{\sf W}$ в исходных условиях, а также $U:=\ln |f| $ для утверждения \eqref{fgiiC}  теоремы \ref{thMR+}. Тогда эквивалентность утверждений    \eqref{fgiC}--\eqref{fgiiiC}  теоремы \ref{thMR+} --- это эквивалентности \eqref{fgiCs}--\eqref{fgiiiCs}
теоремы  \ref{thMR+m}, а дополнение для последовательности ${\sf W}$,  отделённой  углами от мнимой оси, с утверждением (\ref{fgiCS}$'$) соответствует  дополнению с утверждением (\ref{fgiCS+}$'$) из теоремы  \ref{thMR+m} для меры $\upmu$, отделённой углами от мнимой оси.
\end{proof}

Субгармоническая версия теорем {\rm Kh} и \ref{thKh2}  --- следующая  
\begin{theorem}\label{Kh1+S}
Если мера $\upnu\in \Meas^+(\CC)$ отделена углами от  $i\RR$, а  $M\in \sbh_*(\CC)$  конечного типа  $\type_1[M]<+\infty$, то 
 эквивалентны   утверждения:

\begin{enumerate}[{\rm (i)}]
\item\label{khis} для любой функции $u\in \sbh_*(\CC)$ с мерой Рисса 
$\frac{1}{2\pi}\Delta u=\upnu$ и  для любого числа $P\in \RR^+$
найдутся целая функция  $h\in \Hol_*(\CC)$
и множество  $E\in \mathcal B( \RR)$, для которых $\type_1\bigl[u+\ln|h|\bigr]<+\infty$ и выполнено \eqref{ME}; 
 \item\label{khiis} существует  строго возрастающая  последовательность чисел $\{r_n\}_{n\in \NN}\subset \RR^+$, для которой 
$\lim\limits_{n\to\infty}{r_{n+1}}/{r_n}<+\infty$ и   
\begin{equation}\label{lJs}
 \limsup_{N\to  \infty}\sup\limits_{n\leq N}
\Bigl(l_{\upnu}(r_n,r_N)-J_{i\RR}(r_n,r_N; M)\Bigr)<+\infty.
\end{equation}
\end{enumerate}
\end{theorem}
\begin{proof} Пусть $\upmu:=\frac{1}{2\pi}\Delta M \in \Meas^+(\CC)$ --- мера Рисса функции $M\in \sbh_*(\CC)$. По теореме Вейерштрасса\,--\,Адамара\,--\,Линд\-е\-л\-ё\-фа\,--\,Брело мера Рисса $\upmu$ удовлетворяет условию Линделёфа рода $1$ и, тем более, условию \eqref{rlR}[{\rm R}] теоремы \ref{thMR+m}. Эквивалентность 
\eqref{khis}$\Leftrightarrow$\eqref{khiis} теоремы \ref{Kh1+S}следует из  эквивалентности \eqref{fgiCs}$\Leftrightarrow$\eqref{fgivCs} теоремы \ref{thMR+m} и предложения \ref{lemJl} с \eqref{{Jll}m}.
\end{proof}

\begin{proof}[Доказательство теоремы {\rm \ref{thKh2}}]
Положим $\upnu:=n_{\sf Z}$,  $M:=\ln |g|$  и ${\sf W}:=\Zero_g$ в исходных условиях, где ${\sf W}$ удовлетворяет условию Линделёфа рода $1$. Если выполнено \eqref{lJ}, то имеем \eqref{lJs} даже для любой последовательности 
 $\{r_n\}_{n\in \NN}$ со свойствами из утверждения \eqref{khiis} теоремы \ref{Kh1+S}. Пусть $f_{\sf Z}\in \Hol_*(\CC)$ --- какая нибудь целая функция с последовательностью нулей ${\Zero}_f={\sf Z}$. 
Из импликации \eqref{khiis}$\Rightarrow$\eqref{khis}  теоремы \ref{Kh1+S}
для субгармонической функции $u:=\ln |f_{\sf Z}|$ с мерой Рисса $\frac{1}{2\pi}\Delta u=n_{\sf Z}$ для любого числа $P\in \RR^+$
найдутся целая функция  $h\in \Hol_*(\CC)$
и множество  $E\in \mathcal B( \RR)$, для которых выполнено  \eqref{ME}. Если положим $f:=f_{\sf Z}h$, то получаем утверждение
\eqref{fgiCSKh}  теоремы   \ref{thKh2} с  \eqref{fEP}. Обратно, пусть выполнено 
утверждение \eqref{fgiCSKh}  теоремы   \ref{thKh2}. Для 
 любой ц.ф.э.т.  $G$, обращающейся в нуль на ${\sf W}$, функция $h:=G/g$ также ц.ф.э.т. и для пары ц.ф.э.т. $G$ и  $F=fh$ выполнено соотношение  \eqref{fgiC} теоремы \ref{thMR+} с $G$ и $F$ вместо $g$ и $f$ соответственно. 
Из импликации \eqref{fgiC}$\Rightarrow$\eqref{fgiiiC}  теоремы \ref{thMR+},
получаем соотношение \eqref{ZldC}. Из него по соотношению \eqref{{Jll}m} предложения \ref{lemJl} сразу следует 
утверждение \eqref{khii} теоремы  {\rm Kh} с  \eqref{lJ}.
\end{proof}

\section{Характеристики последовательностей, мер и зарядов}
\setcounter{equation}{0}
\subsection{Свойства логарифмической функции интервалов для  зарядов}\label{2_3_1} 
Определим некоторые считающие радиальные функции с весом. В  \cite[теорема A]{Kha99} они рассматривались в значительно более общей форме для последовательностей точек и множеств в $\CC^n$.  В связи с тем, что в нашей статье рассматривается рост функций вдоль $i\RR$, для точек  $z\in \CC$ \textit{главные значения аргументов\/} $\arg z$  удобнее выбирать из интервала $[-\pi/2, 3\pi /2)$.

\begin{definition}[{\cite[определение 3]{Kha20}}]
Пусть $\upnu\in \Meas(\CC)$, $k\colon \RR \to \RR$  ---  {\it $2\pi$-периодическая борелевская   функция.\/} Введем в рассмотрение \textit{считающую  функцию заряда\/  $\upnu$ с весом\/ $k$} как интеграл 
 \begin{equation}\label{c)}
\upnu(r;k):=\int_{\overline D(r)} k (\arg z)\dd \upnu(z), \quad r\in \RR^+.
\end{equation}
При $k\equiv 1$, очевидно, $\upnu(r;1)\overset{\eqref{df:nup}}{\equiv} \upnu^{\rad}(r)$ для $r\in \RR^+$.
\end{definition}

В частных случаях $k=\cos^{\pm}$, т.\,е.
при $\cos^+:=\sup\{\cos ,0\}$, $\cos^-:=\sup\{-\cos ,0\}$, 
 из определений \eqref{df:lm} в обозначении \eqref{c)} 
интегрированием по частям  при  $0< r < R <+\infty$   получаем 
\begin{subequations}\label{l_mu}
\begin{align}
l_{\upnu}^{\rght}(r, R)&\overset{\eqref{df:dDlm+}}{=}
\int_r^R \frac{\dd \upnu(t;\cos^+)}{t}
\notag\\
&=\frac{\upnu (R;\cos^+)}{R}-\frac{\upnu (r;\cos^+)}{r}
+\int_r^R\frac{\upnu (t;\cos^+)}{t^2} \dd t, 
\tag{\ref{l_mu}r}\label{l_mu_m+}
\\
l_{\upnu}^{\lft}(r, R)&\overset{\eqref{df:dDlm-}}{=}
\int_r^R \frac{1}{t} \dd \upnu(t;\cos^-)
\notag \\
&=\frac{\upnu (R;\cos^-)}{R}-\frac{\upnu (r;\cos^-)}{r}
+\int_r^R\frac{\upnu (t;\cos^-)}{t^2} \dd t.  
\tag{\ref{l_mu}l}\label{l_mu_m-}
\end{align}
\end{subequations}
Положим 
\begin{subequations}\label{l_mu-}
\begin{align}
\breve l_{\upnu}^{\rght}(r, R)&\overset{\eqref{l_mu_m+}}{:=}\int_r^R\frac{\upnu (t;\cos^+)}{t^2} \dd t, \quad 0< r < R < +\infty,
\tag{\ref{l_mu-}r}\label{l_mu_m-+}
\\
\breve l_{\upnu}^{\lft}(r, R)&\overset{\eqref{l_mu_m-}}{:=}\int_r^R\frac{\upnu (t;\cos^-)}{t^2} \dd t,\quad 0< r < R < +\infty,  
\tag{\ref{l_mu-}l}\label{l_mu_m--}
\\
\intertext{а для \textit{положительной меры\/} $\upmu\in \Meas^+(\CC)$ ещё и}
\breve l_{\upmu}(r, R)&\overset{\eqref{df:dDlLm}}{:=}\max \{ \breve l_{\upmu}^{\lft}(r, R), \breve l_{\upmu}^{\rght}(r, R)\}, \quad 0< r < R < +\infty .
\tag{\ref{l_mu-}m}\label{ml}
\end{align}
\end{subequations}
В частности, для произвольной последовательности ${\sf Z}\subset \CC$ со считающей мерой $n_{\sf Z}$ из \eqref{df:divmn} полагаем 
\begin{equation}\label{lZ}
\breve l_{\sf Z}^{\rght}(r, R)\overset{\eqref{l_mu_m-+}}{:=}\breve l_{n_{\sf Z}}^{\rght}(r, R), \;
\breve l_{\sf Z}^{\lft}(r, R)\overset{\eqref{l_mu_m--}}{:=}\breve l_{n_{\sf Z}}^{\lft}(r, R), \;\breve l_{\sf Z}(r,R)\overset{\eqref{ml}}{:=}\breve l_{n_{\sf Z}}(r,R).
\end{equation}

Для заряда $\upnu\in \Meas(\CC)$ \textit{центрально симметричный\/} ему  \begin{equation}\label{nu-}
\begin{split}
\text{\it заряд } &\upnu_{\times}\in \Meas(\CC) \quad\text{\it определяется равенствами}\\ 
&\upnu_\times(B):=\upnu (-B) \quad\text{\it по всем } B\in \mathcal B_{\rm b}(\CC),
\end{split}
\end{equation}
а \textit{зеркально симметричный\/} ему {\it относительно $i\RR$\/}  
\begin{equation}\label{nus}
\begin{split}
\text{\it заряд } &\upnu_{\leftrightarrow}\in \Meas(\CC) \quad\text{\it определяется равенствами}\\ 
&\upnu_\leftrightarrow(B):=\upnu (-\bar B) \quad\text{\it по всем } B\in \mathcal B_{\rm b}(\CC).
\end{split}
\end{equation}
Очевидно, операции $\upnu\mapsto \upnu_\times$ и $\upnu\mapsto \upnu_\leftrightarrow$ как любые симметрии  --- линейные инволюции. 
Заряд $\upnu\in \Meas(\CC)$  \textit{чётный,\/} если 
$\upnu_\times=\upnu$; {\it нечётный,\/} если $\upnu_\times=-\upnu$; $i\RR$-\textit{зеркальный,\/} если $\upnu=\upnu_\leftrightarrow$;
$i\RR$-\textit{антизеркальный,\/} если $\upnu=-\upnu_\leftrightarrow$. 
Для любого  заряда $\upnu\in \Meas(\CC) $ заряд $\upnu+\upnu_\times$ чётный,  
$\upnu-\upnu_\times$ нечётный, а  $\upnu+\upnu_\leftrightarrow$ ---  $i\RR$-зеркальный, 
$\upnu-\upnu_\leftrightarrow$ ---  $i\RR$-антизеркальный. Таким образом,
имеют место стандартные представления заряда $\upnu$ в виде сумм 
\begin{equation}\label{reprnu}
\upnu=\frac12(\upnu+\upnu_{\times})+\frac12(\upnu-\upnu_{\times}),
\quad \upnu=\frac12(\upnu+\upnu_{\leftrightarrow})+\frac12(\upnu-\upnu_{\leftrightarrow})
\end{equation} 
чётного $\frac12(\upnu+\upnu_{\times})$ и нечётного $\frac12(\upnu-\upnu_{\times})$, а также  $i\RR$-зеркального $\frac12(\upnu+\upnu_{\leftrightarrow})$ и $i\RR$-антизеркального $\frac12(\upnu-\upnu_{\leftrightarrow})$ зарядов соответственно. 

Из определений логарифмических функций/(суб)мер интервалов имеем 

\begin{propos}\label{prcc} Для любого заряда $\upnu\in \Meas(\CC)$ при всех значениях 
$0<r<R<+\infty$ в обозначениях \eqref{nu-}--\eqref{nus}  имеют место тождества
\begin{subequations}\label{vnu}
\begin{align}
l_{\upnu}^{\rght}(r, R)&\equiv l_{\upnu_\times}^{\lft}(r, R)
\equiv l_{\upnu_\leftrightarrow}^{\lft}(r, R),
\tag{\ref{vnu}r}\label{{vnu}r}\\
l_{\upnu}^{\lft}(r, R)&\equiv l_{\upnu_\times}^{\rght}(r, R)
\equiv l_{\upnu_\leftrightarrow}^{\rght}(r, R),
\tag{\ref{vnu}l}\label{{vnu}l}\\
\intertext{а для положительной меры  $\upmu\in \Meas^+(\CC)$ также тождества}
l_{\upmu}(r, R)&\equiv l_{\upmu_\times}(r, R)
\equiv l_{\upmu_\leftrightarrow}(r, R).
\tag{\ref{vnu}m}\label{{vnu}m}
\end{align}
\end{subequations}
Эти тождества остаются в силе, если всюду в них  логарифмические функции/(суб)меры интервалов  $l_{\cdot}^{\cdot}(\dots)$  одновременно заменить  на соответствующие им функции вида $\breve l_{\cdot}^{\cdot}(\dots)$ из  \eqref{l_mu-}--\eqref{lZ}. 
\end{propos}

\begin{propos}[{\cite[предложение 3]{Kha20}}]\label{pr:lm}
Пусть $\upnu\in \Meas(\CC)$ и $\upmu\in \Meas^+(\CC)$ конечной верхней плотности при порядке $1$, $0<r_0\in \RR^+$. Тогда 
\begin{equation}\label{ml+-}
\begin{cases}
\bigl|\breve l_{\upnu}^{\rght}(r, R)-l_{\upnu}^{\rght}(r, R)\bigr|\overset{\eqref{l_mu_m-+}}{=}O(1)\\
\bigl|\breve l_{\upnu}^{\lft}(r, R)-l_{\upnu}^{\lft}(r, R)\bigr|
\overset{\eqref{l_mu_m--}}{=}
O(1) \\
\bigl|\breve l_{\upmu}(r,R) - l_{\upmu}(r,R)\bigr |\overset{\eqref{ml}}{=}O(1) 
\end{cases}
\text{для всех  $r_0\leq r<R<+\infty$,}
\end{equation}
 для  любых фиксированных чисел  $a\in (0,1]$, $b\in [1,+\infty)$  
\begin{equation}\label{ml+ab}
\begin{cases}
\bigl| l_{\upnu}^{\rght}(r, R)-l_{\upnu}^{\rght}(ar,bR)\bigr|=O(1)\\
\bigl| l_{\upnu}^{\lft}(r, R)-l_{\upnu}^{\lft}(ar, bR)\bigr|=O(1) \\
\bigl|l_{\upmu}(r,R) - l_{\upmu}(ar,bR)\bigr |=O(1)
\end{cases}
\text{для всех  $r_0\leq r<R<+\infty$,}
\end{equation}
а при соглашениях 
\begin{equation}\label{pml}
l_{\upnu}^{\rght}(1,r_0):=-l_\upnu^{\rght}(r_0,1), \;
l_{\upnu}^{\lft}(1,r_0):=-l_\upnu^{\lft}(r_0,1), \;
l_{\upnu}(1,r_0):=-l_\upnu(r_0,1) 
\end{equation}
всегда имеем 
\begin{equation}\label{1r0}
\bigl|l_{\upnu}^{\rght}(1,r_0)\bigr|+
\bigl|l_{\upnu}^{\lft}(1,r_0)\bigr|+
\bigl|l_{\upmu}(1,r_0)\bigr|<+\infty.
\end{equation}
\end{propos}

\begin{remark}\label{rembr}
Из соотношений \eqref{ml+-} предложения \ref{pr:lm} всюду ниже в определениях и утверждениях с зарядом,  мерой, последовательностью точек  \textit{при условии конечной  верхней плотности\/} для них  различные {\it логарифмические функции\/}  или {\it (суб)меры интервалов\/}
вида $l_{\cdot}^{\cdot}(\dots)$ из  \eqref{df:lm} и \eqref{df:l}  можно заменить  на соответствующие им \textit{логарифмические функции/(суб)меры  интервалов\/} вида $\breve l_{\cdot}^{\cdot}(\dots)$ из  \eqref{l_mu-}--\eqref{lZ}, которые в дальнейшем 
обозначаем так же через $l_{\cdot}^{\cdot}(\dots)$ без верхнего  математического акцента $\,\breve{\phantom{l}}\,$. Кроме того, при использовании логарифмических   функций/(суб)мер интервалов выполнение с фиксированным $r_0\in \RR^+\setminus \{0\}$ каких-либо соотношений при всех $r_0\leq r<R<+\infty$ с точностью до аддитивного слагаемого $C\in \RR$ в силу \eqref{1r0} при соглашении \eqref{pml} можно заменить на их выполнение  при всех $1\leq r<R<+\infty$, т.\,е. положить  $r_0:=1$.
\end{remark}

Из соотношений \eqref{ml+ab} предложения \ref{pr:lm} и  замечания \ref{rembr} легко следует

\begin{propos}\label{pr:rn}
Пусть $\upnu\in \Meas^+(\CC)$ и\/  $\upmu\in \Meas^+(\CC)$ конечной верхней плотности при порядке $1$. Отношение
$\upnu\overset{\eqref{lJMmul}}{\curlyeqprec}_{\RR} \upmu$ имеет место тогда и только тогда, когда выполнено утверждение \eqref{fgivCs}
теоремы\/ {\rm \ref{thMR+m}} с \eqref{cprec}.
\end{propos}

\subsection{Различные формы условия Бляшке для зарядов}\label{cBl}

В \cite[следствие 4.1, теорема 4]{KhI} показано, что классическое выметание \textit{положительной меры\/}  $\upmu \in \Meas^+ (\CC)$ из  полуплоскости или углов  возможно тогда и только тогда, когда для меры $\upmu$
выполнено {\it классическое условие Бляшке\/} \cite[4.1, определения 4.1, 4.6]{KhI},  которое в {\it правой\/} $\CC_{\rght}$ 
или   {\it левой полуплоскости\/} $\CC_{\lft}$ 
для \textit{ заряда\/} $\upnu \in \Meas(\CC)$ означает, что соответственно
\begin{equation}\label{Blcl}
\int_{\CC_{\rght}\setminus \DD} \Re \frac{1}{z}\dd |\upnu| (z)<+\infty
\quad\text{или}\quad
 \int_{\CC_{\lft}\setminus \DD} \Bigl(-\Re \frac{1}{z}\Bigr)\dd |\upnu| (z)<+\infty
\end{equation}
где подынтегральные выражения \textit{положительны.\/}
В терминах правых и левых логарифмических функций интервалов условия \eqref{Blcl}  эквивалентны соответственно  соотношениям
\begin{equation}\label{lchar}
l_{|\upnu|}^{\rght}(1,+\infty)\overset{\eqref{df:dDlLm}}{<}+\infty\quad\text{или}\quad l_{|\upnu|}^{\lft}(1,+\infty)<+\infty. 
\end{equation}
Заряд $\upnu\in \Meas(\CC)$ удовлетворяет {\it условию Бляшке вне\/} $i\RR$, если 
\begin{equation}\label{nuB}
\Bigl(l_{|\upnu|}(1,+\infty)<+\infty\Bigr)\overset{\text{или}}{\Longleftrightarrow}\left(\int_{\CC\setminus \DD}\left|\Re \frac{1}{z}\right| \dd |\upnu|(z)<+\infty\right).
\end{equation} 
В частности, последовательность ${\sf Z}=\{{\sf z}_k\}\subset \CC$
удовлетворяет  условию Бляшке вне $i\RR$, если её считающая мера $n_{\sf Z}$ 
удовлетворяет \eqref{nuB}, т.е.
\begin{equation*}
\Bigl(l_{\sf Z}(1,+\infty)\overset{\eqref{df:dDlL}}{<}+\infty\Bigr)
\overset{\text{или}}{\Longleftrightarrow}\left(\sum_{k} \left|\Re \frac{1}{{\sf z}_k}\right| <+\infty\right).
\end{equation*}

Для заряда $\upnu$ на $\CC$ конечной верхней плотности очевидно 
\begin{equation}\label{rint}
l_{|\upnu|}(1,r):=\max\left\{l_{|\upnu|}^{\rght}(1,r), 
 l_{|\upnu|}^{\lft}(1,r)\right\}\leq \int_{1}^{r} \frac{|\upnu|^{\rad}(t)}{t^2}\dd t
+O(1)
\end{equation}
при $ r\to +\infty$,  откуда сразу следует 
\begin{propos}\label{prcls}  Если заряд $\upnu\in \Meas(\CC)$
из класса сходимости около $\infty$ при порядке $1$ \cite[2.1]{KhI}, т.\,е. 
при  $r=+\infty$ конечен интеграл из правой части \eqref{rint}, то заряд $\upnu$ удовлетворяет  условию Бляшке \eqref{nuB} вне $i\RR$.
\end{propos}
Выметание конечного рода $q\in \NN$ \textit{заряда\/} $\upnu\in \Meas(\CC)$ из полуплоскости или угла обеспечивалось в 
\cite[теоремы A1, 3, 5, следствие 3.1]{KhII} при $q=1$ {\it условием Бляшке рода\/ $1$ для заряда\/} $\upnu$ \cite[(1.5), (3.23), (3.42), {\bf 4.3}]{KhII}, которое в  {\it полуплоскостях\/}  $\CC_{\rght}$ 
или $\CC_{\lft}$  означает, что соответственно 
\begin{equation}\label{Blclq}
\sup_{r\geq 1} \left|\int_{S(r)}\Re \frac{1}{z}\dd \upnu (z)\right|<+\infty
\quad\text{для}\quad S(r):=\overline D(r)\cap \left[
  \begin{array}{c}
 \CC_{\rght}\setminus \DD \\
\CC_{\lft}\setminus \DD\\
  \end{array}
\right. .
\end{equation}
В терминах правых и левых логарифмических функций интервалов условия \eqref{Blclq} эквивалентны соответственно  соотношениям
\begin{equation}\label{Bl}
\sup_{r\geq 1} \left|l_{\upnu}^{\rght}(1,r)\right|<+\infty \quad\text{или}\quad
\quad \sup_{r\geq 1} \left|l_{\upnu}^{\lft}(1,r)\right|<+\infty.
\end{equation}
В случае  одновременного выметания рода $1$ заряда $\upnu$ на мнимую ось
из  $\CC_{\rght}$ и $\CC_{\lft}$ можно ещё более ослабить условия \eqref{Blclq}--\eqref{Bl}.

\begin{definition}\label{df:bB}
Заряд $\upnu\in \Meas(\CC)$ 
удовлетворяет {\it двустороннему условию Бляшке\footnote{В \cite[определение 3.3]{kh91AA}, \cite[определение 2.2.1]{KhDD92}  оно называлось условием A.} рода\/ $1$ вне  мнимой оси\/ $i\RR$,\/} 
если 
\begin{equation}\label{cB2}
\sup_{r\geq 1} \left| \int_{\overline D(r)\setminus \DD}\Bigl|\Re\frac{1}{z} \Bigr| \dd \upnu(z)\right|< +\infty,
\end{equation}
что в терминах логарифмических функций интервалов записывается  как 
\begin{equation}\label{cB2l}
\sup_{r\geq 1}\bigl|l_{\upnu}^{\rght}(1,r)+l_{\upnu}^{\lft}(1,r)\bigr|<+\infty.
\end{equation}
Для любого нечётного или $i\RR$-антизеркального  заряда $\upnu\in \Meas(\CC)$ в \eqref{cB2}--\eqref{cB2l} выражения под модулем после $\sup_{r\geq 1}$ равны нулю. Следовательно,  для таких зарядов   выполнено  двустороннее условие Бляшке.
\end{definition}

\begin{propos}\label{pr2}  Имеют место следующие взаимосвязи между различными формами условия Бляшке: 
\begin{enumerate}[{\rm [{b}1]}]
\item\label{B1} Из классического условия Бляшке \eqref{Blcl}--\eqref{lchar} в 
правой или  в левой полуплоскости следует условие Бляшке \eqref{Blclq}--\eqref{Bl} рода $1$ соответственно в правой или левой полуплоскости. Кроме того, 
\begin{enumerate}[{\rm (i)}]
\item\label{B1i}  для  положительных мер верно и обратное; 
\item\label{B1ii}  существуют заряды $\upnu$ нулевого типа $\type_1[\upnu]\overset{\eqref{fden}}{=}0$, и  обязательно  из класса расходимости около $\infty$ (при порядке $1$), т.\,е. с расходящимся при $r=+\infty$  интегралом 
из правой части \eqref{rint}, удовлетворяющие условию Бляшке \eqref{Blclq}--\eqref{Bl}  рода\/ $1$ одновременно и в $\CC_{\rght}$, и в $\CC_{\lft}$, но не удовлетворяющей классическим условиям  Бляшке \eqref{Blcl}--\eqref{lchar} ни в $\CC_{\rght}$, ни в $\CC_{\lft}$.
\end{enumerate}
\item\label{B2} Из условия Бляшке \eqref{Blclq}--\eqref{Bl} рода\/ $1$ одновременно в $\CC_{\rght}$ и в $\CC_{\lft}$ следует двустороннее условие Бляшке \eqref{cB2}--\eqref{cB2l}. 
\begin{enumerate}[{\rm (i)}]
\item\label{B2i}  для  положительных мер верно и обратное; 
\item\label{B2ii}  существуют заряды нулевого типа, и  обязательно  из класса расходимости около $\infty$, удовлетворяющие двустороннему усло\-в\-ию Бляшке \eqref{cB2}--\eqref{cB2l}, но не удовлетворяющие  условиям  Бля\-шке \eqref{Blclq}--\eqref{Bl} рода $1$ ни в $\CC_{\rght}$, ни в $\CC_{\lft}$;
\end{enumerate}
\item\label{B3} Для $\upnu\in \Meas(\CC)$ 
следующие три утверждения эквивалентны:
\begin{enumerate}[{\rm (i)}]
\item\label{B3i} $\upnu$ удовлетворяет двустороннему условию Бляшке \eqref{cB2}--\eqref{cB2l}; 
\item\label{B3ii} чётный  заряд $\upnu+\upnu_\times$, где заряд $\upnu_\times$, центрально симметричный заряду $\upnu$, определён в \eqref{nu-},  удовлетворяет условию Бля\-ш\-ке рода\/ $1$ как в $\CC_{\rght}$, так и в $\CC_{\lft}$ в смысле \eqref{Blclq}--\eqref{Bl};
\item\label{B3iii}  заряд $\upnu+\upnu_\leftrightarrow$, где заряд $\upnu_\leftrightarrow$, зеркально симметричный заряду $\upnu$ относительно $i\RR$, определён в \eqref{nus},  удовлетворяет условию Бляшке рода\/ $1$ как в $\CC_{\rght}$, так и в $\CC_{\lft}$ в смысле \eqref{Blclq}--\eqref{Bl}.
\end{enumerate}
\end{enumerate}
\end{propos}
\begin{proof} Первое утверждение в {b}\ref{B1} следует из неравенств 
\begin{equation}\label{B1in}
\bigl|l_{\upnu}^{\rght}(r_0,r)\bigr|\leq l_{|\upnu|}^{\rght}(r_0,r), 
\quad \bigl|l_{\upnu}^{\lft}(r_0,r)\bigr|\leq l_{|\upnu|}^{\lft}(r_0,r) 
\end{equation}
и определений  в форме \eqref{lchar} и \eqref{Bl}. Для положительной меры 
$\upnu=|\upnu|$ в \eqref{B1in} имеем равенства, что даёт утверждение 
{b}\ref{B1i}. 

Необходимые требования к построению заряда $\upnu$ из {b}\ref{B1ii} о принадлежности к классу расходимости --- следствие предложения \ref{prcls}. Для построения заряда $\upnu$ с требуемыми в {b}\ref{B1ii} свойствами рассмотрим 
\begin{example}\label{ex1} Пусть  $m\colon \RR^+\to \RR^+$ --- возрастающая функция нулевого типа  около $+\infty$  из класса расходимости при порядке $1$ \cite[2.1]{KhI}. К примеру, можно выбрать   
$m(x)\equiv x/\ln(e+x)$, $x\in \RR^+$. Рассмотрим положительную меру $\upmu$ с  $\supp \upmu \subset \RR$ и с функцией распределения $\upmu^{\RR}(x)\overset{\eqref{nuR}}{\equiv} m(x)$ при $x\in \RR^+$ и $\upmu^{\RR}(x)\equiv -m(-x)$ при $x\in \RR^+$. Кроме того, рассмотрим <<поворот>> меры $\upmu$ на угол $\theta \in (0,\pi/2)$, обозначаемый как $\upmu_{\theta}$ и определённый по правилу 
$\upmu_{\theta}(S):=\upmu(e^{-i\theta}S)$ для всех $S\in \mathcal B_{\rm b}(\CC)$.
Тогда для заряда $\upnu:=\upmu-\upmu_{\theta}$ имеем 
\begin{equation}\label{lnumu}
\begin{split}
&l_{\upnu}^{\rght}(1,r)\equiv l_{\upnu}^{\lft}(1,r) \equiv 0 \quad
\text{при всех $1<r<+\infty$},\\ 
&l_{|\upnu|}^{\rght}(1,r)\equiv l_{|\upnu|}^{\lft}(1,r) \equiv 
l_{2\upmu}^{\rght}(1,r)\to +\infty\quad\text{при  }r\to +\infty, 
\end{split}
\end{equation}
поскольку мера $\upmu$ из класса расходимости \cite[предложение 2.1]{KhI}. 
\end{example}
Пример  \ref{ex1} по определениям  \eqref{Bl} и \eqref{lchar} в силу соотношений \eqref{lnumu} даёт  требуемый в  {b}\ref{B1ii} заряд $\upnu$.

Первое утверждение в {b}\ref{B2} следует из очевидного  неравенства 
\begin{equation}\label{2cB}
\bigl|l_{\upnu}^{\rght}(1,r)+l_{\upnu}^{\lft}(1,r)\bigr|
\leq \bigl|l_{\upnu}^{\rght}(1,r)\bigr|+\bigl|l_{\upnu}^{\lft}(1,r)\bigr|
\text{ для $1<r<+\infty$}
\end{equation}
и определений в форме \eqref{Bl} и \eqref{cB2l}.  Для положительной меры 
$\upnu=|\upnu|$ в \eqref{2cB} имеем равенства, что даёт утверждение 
{b}\ref{B2i}. 

Для {b}\ref{B2ii}   достаточно рассмотреть любой нечётный заряд $\upnu$ с 
положительным сужением на правую полуплоскость, принадлежащим к классу расходимости. Например, можно взять нечётный заряд $\upnu$, сужение которого на $\CC_{\rght}$ совпадает с сужением меры  $\upmu$ из примера \ref{ex1} на $\CC_{\rght}$.

Эквивалентность трёх утверждений п. {b}\ref{B3}
следует из определений  в формах  \eqref{cB2l} и \eqref{Bl}, а также из тождеств \eqref{vnu} предложения \ref{prcc}.
\end{proof}

\subsection{Условие Линделёфа для зарядов, мер и комплексных последовательностей}\label{cL}
Очевидно, чётный  заряд $\upnu$ удовлетворяет  условию Линделёфа рода $1$.  Из условия Линделёфа \eqref{con:Lp+}, вообще говоря, не следует ни одно из   условий Бляшке, упоминавшихся  в \S~\ref{cBl}, как и наоборот, --- 
ни одно из  этих  условий 
не влечёт за собой условие Линделёфа \eqref{con:Lp+}. Но 
условие Бляшке \eqref{Blclq}--\eqref{Bl} рода\/ $1$ одновременно в $\CC_{\rght}$ и в $\CC_{\lft}$ при 
\begin{equation}\label{con:LpIm} 
\sup_{r\geq 1}\Biggl|\int_{\overline D(r)\setminus \overline \DD} \Im \frac{1}{z}
\dd \upnu(z) \Biggr|<+\infty
\end{equation} 
даёт  условие Линделёфа \eqref{con:Lp+} для заряда $\upnu\in \Meas(\CC)$. 
Для  последовательностей ${\sf Z}$ из \eqref{Z} условие 
  \eqref{con:LpIm}  для считающей меры $n_{\sf Z}$ записывается как
\begin{equation}\label{con:LpIZ} 
\sup_{r\geq 1}\Biggl|\sum_{1< |{\sf z}_k|\leq r} \Im \frac{1}{{\sf z}_k} \Biggr|<+\infty,
\end{equation}

Для  заряда $\upnu\in \Meas(\CC)$ и числа $\alpha \in \RR$
его {\it повороты против часовой стрелки $\upnu_{\circlearrowleft \alpha}$} и 
{\it  по часовой стрелке $\upnu_{\circlearrowright \alpha}$ на угол $\alpha$} определим как  заряды  
\begin{equation}\label{pnZ}
\upnu_{\circlearrowleft \alpha}(B):=\upnu(e^{-i\alpha}B), \quad 
\quad \upnu_{\circlearrowright \alpha}(B):=\upnu_{\circlearrowleft (-\alpha)}(B)
\quad \text{для $B\in {\mathcal B}_{\rm b}(\CC)$}. 
\end{equation}
В частности, $\upnu_{\times}\overset{\eqref{nu-}}{=}\upnu_{\circlearrowleft \pi}
=\upnu_{\circlearrowright \pi}$. Свойство \eqref{con:LpIm} 
можно записать в виде 
$\sup_{r\geq 1}\bigl|l_{\upnu_{\circlearrowleft \pi/2}}^{\rght}(1,r)-l_{\upnu_{\circlearrowleft \pi/2}}^{\lft}(1,r)\bigr|<+\infty$,
а условие Линделёфа \eqref{con:Lp+} рода $1$ эквивалентно паре соотношений 
\begin{equation}\label{LimL}
\begin{cases}
\sup_{r\geq 1}\bigl|l_{\upnu}^{\rght}(1,r)-l_{\upnu}^{\lft}(1,r)\bigr|<+\infty,\\
\sup_{r\geq 1}\bigl|l_{\upnu_{\circlearrowleft \pi/2}}^{\rght}(1,r)-l_{\upnu_{\circlearrowleft \pi/2}}^{\lft}(1,r)\bigr|<+\infty,
\end{cases}
\end{equation}
Используя  эту эквивалентность с \eqref{LimL}, для мер $\upmu\in \Meas^+(\CC)$ получаем 
\begin{propos}\label{remL}
Мера $\upmu\in \Meas^+(\CC)$ удовлетворяет условию Линделёфа   
\eqref{con:Lp+} рода $1$  тогда и только тогда, когда
одновременно 
\begin{equation}\label{lL}
\begin{split}
\sup_{1\leq r<R<+\infty} 
\Bigl\{\bigl|l_{\upmu}^{\rght}(r,R)-l_{\upmu}^{\lft}(r,R)\bigr|,& \; \bigl|l_{\upmu}^{\rght}(r,R)-l_{\upmu}(r,R)\bigr|\Bigr\}<+\infty,
\\
\sup_{1\leq r<R<+\infty} 
\Bigl\{\bigl|l_{\upmu_{\circlearrowleft \pi/2}}^{\rght}(r,R)-l_{\upmu_{\circlearrowleft \pi/2}}^{\lft}&(r,R)\bigr|,\Bigr.\\\Bigl.
 \bigl|l_{\upmu_{\circlearrowleft \pi/2}}^{\rght}&(r,R)-l_{\upmu_{\circlearrowleft \pi/2}}(r,R)\bigr|\Bigr\}<+\infty.
\end{split}
\end{equation}
\end{propos}

 Для комплексной последовательности ${\sf Z}$  поворотам  \eqref{pnZ} её считающей меры 
$n_{\sf Z}$ соответствуют {\it повороты последовательности ${\sf Z}$ против часовой стрелки ${\sf Z}_{\circlearrowleft \alpha}$} и 
{\it  по часовой стрелке ${\sf Z}_{\circlearrowright \alpha}$ на угол $\alpha$}, которые  по п.~\ref{1_2_4}  можно определить  соответственно как  последовательности 
${\sf Z}_{\circlearrowleft \alpha}=e^{i\alpha}{\sf Z}$ и  
${\sf Z}_{\circlearrowright \alpha}=e^{-i\alpha}{\sf Z}$. В частности, 
${\sf Z}_{\circlearrowleft \pi/2}=i{\sf Z}$ и  
${\sf Z}_{\circlearrowright \pi/2}=-i{\sf Z}$ и  свойство \eqref{con:LpIZ} можно записать как 
$\sup_{r\geq 1}\bigl|l_{i{\sf Z}}^{\rght}(1,r)-l_{i{\sf Z}}^{\lft}(1,r)\bigr|<+\infty$,
а условие Линделёфа \eqref{con:LpZ}  рода $1$ эквивалентно паре соотношений 
\begin{equation*}
\begin{cases}
\sup_{r\geq 1}\bigl|l_{\sf Z}^{\rght}(1,r)-l_{\sf Z}^{\lft}(1,r)\bigr|<+\infty,\\
\sup_{r\geq 1}\bigl|l_{i{\sf Z}}^{\rght}(1,r)-l_{i{\sf Z}}^{\lft}(1,r)\bigr|<+\infty.
\end{cases}
\end{equation*}
Предложение \ref{remL} при этом даёт 
\begin{propos}\label{remLZ}
Последовательность  ${\sf Z}\subset \CC$ удовлетворяет условию Линделёфа   
\eqref{con:Lp+} рода $1$  тогда и только тогда, когда
одновременно 
\begin{equation*}
\begin{split}
\sup_{1\leq r<R<+\infty} 
\Bigl\{\bigl|l_{\sf Z}^{\rght}(r,R)-l_{\sf Z}^{\lft}(r,R)\bigr|,& \;\bigl|l_{\sf Z}^{\rght}(r,R)-l_{\sf Z}(r,R)\bigr|\Bigr\}<+\infty,
\\
\sup_{1\leq r<R<+\infty} 
\Bigl\{\bigl|l_{i\sf Z}^{\rght}(r,R)-l_{i\sf Z}^{\lft}(r,R)\bigr|,&\;
 \bigl|l_{i\sf Z}^{\rght}(r,R)-l_{i\sf Z}(r,R)\bigr|\Bigr\}<+\infty.
\end{split}
\end{equation*}
\end{propos}

\begin{propos}\label{prL} Пусть мера $\upmu\in \Meas^+ (\CC)$ конечной верхней плотности при порядке $1$ и выполнено условие \eqref{rlR}{\bf [R]} из теоремы\/ {\rm \ref{thMR+m}}. Тогда существует  
мера  $\upbeta \in \Meas^+(i\RR)$, для которой сумма  $\upmu':=\upmu+\upbeta \in \Meas^+(\CC)$ --- мера конечной верхней плотности при порядке $1$, удовлетворяющая условию Линделёфа \eqref{con:Lp+} рода $1$ и, очевидно,   
\begin{equation}\label{mus}
\upmu \curlyeqprec_{|\Re|} \upmu'=\upmu+\upbeta\curlyeqprec_{|\Re|}\upmu.  
\end{equation}
Если $\supp \upmu \subset \CC_{\overline \rght}$ или $\supp \upmu \subset \CC_{\overline \lft}$, то меру $\upmu'\geq \upmu$    с требуемыми свойствами можно задать   и как чётную меру $\upmu'\overset{\eqref{nu-}}{:=}\upmu +\upmu_{\times}$.
\end{propos}
\begin{proof}
В обозначениях из \eqref{pnZ} для поворота меры положим
 \begin{equation}\label{an}
\begin{cases}
b_n^{\rght}:=2^n
\Bigl(
l_{\upmu_{\circlearrowleft \pi/2}}^{\lft}
(2^{n-1},2^n)
- l_{\upmu_{\circlearrowleft \pi/2}}^{\rght}
(2^{n-1},2^n)\Bigr)^+,
\\
b_n^{\lft}:=2^n
\Bigl(l_{\upmu_{\circlearrowleft \pi/2}}^{\rght}
(2^{n-1},2^n)-l_{\upmu_{\circlearrowleft \pi/2}}^{\lft}
(2^{n-1},2^n)\Bigr)^+,
\end{cases}
\quad n\in \NN,
\end{equation}
и для мер Дирака $\updelta_{\pm i2^n}$ в точках $\pm i2^n$ рассмотрим меру 
\begin{equation}\label{alp}
\upbeta :=\sum_{n=1}^{\infty} b_n^{\rght}\updelta_{i2^n} +
\sum_{n=1}^{\infty} b_n^{\lft}\updelta_{-i2^n}\in \Meas^+(i\RR). 
\end{equation}
По построению \eqref{an}--\eqref{alp} и определениям \eqref{df:lm} 
для всех $n\in \NN$ имеем 
\begin{equation}\label{l2n}
l_{\upmu_{\circlearrowleft \pi/2}}(2^{n-1}, 2^n)=l_{(\upmu+\upbeta)_{\circlearrowleft \pi/2}}^{\rght}(2^{n-1}, 2^n)
=l_{(\upmu+\upbeta)_{\circlearrowleft \pi/2}}^{\lft}(2^{n-1}, 2^n),
\end{equation}
а мера $\upbeta\in \Meas^+(\RR)$ конечной верхней плотности при порядке $1$. Из \eqref{l2n} легко следует второе соотношение из \eqref{lL} предложения \ref{remL} для меры $\upmu+\upbeta$, а первое  соотношение из  \eqref{lL} для меры $\upmu+\upbeta$ следует из условия   \eqref{rlR}{\bf [R]}, поскольку $\supp \upbeta \subset i\RR$. Таким образом, по предложению  \ref{remL}  мера $\upmu+\upbeta$ удовлетворяет условию Линделёфа рода $1$. Дополнение для мер  $\supp \upmu \subset \CC_{\overline \rght}$ или $\supp \upmu \subset \CC_{\overline \lft}$ очевидно, поскольку мера $\upmu +\upmu_{\times}$ чётная. 
\end{proof}

\section{Классическое выметание на вещественную и мнимую ось} 
\setcounter{equation}{0}

Пусть $p\in \RR^+$. Напомним, что в \cite[замечание 4.2]{KhI}, \cite[определение 1.4]{Kha91}, \cite[2]{Kha01l} замкнутая cистема  лучей $S$ 
на комплексной с вершиной в нуле называлась {\it $p$-до\-п\-у\-с\-т\-и\-мой\/},
если раствор любого открытого угла, дополнительного к $S$, т.е. связной компоненты в  $\CC\setminus S$,  меньше, чем  $\pi /p$. 
\begin{theoB}[{\rm \cite[основная теорема]{Kha91},  \cite[теорема 1.3.1]{KhDD92}, \cite[предложение 2.1]{Kha01l},  \cite[теорема 8]{KhI}}]\label{BR} Пусть $p\in \RR^+$. Если замкнутая система лучей $S$ на $\CC$ с вершиной в нуле   $p$-допустима, то для любой субгармонической функции $v\in \sbh_*(\CC)$ конечного типа $\type_p[v]<+\infty$ при порядке $p$
существует субгармоническая функция $v^{\bal}$ конечного типа 
$\type_p[v^{\bal}]<+\infty$ при порядке $p$, гармоническая в каждом дополнительном к системе лучей $S$ угле, совпадающая с функцией $v$ на каждом луче из $S$ и при этом $v\leq v^{\bal}$ на всей плоскости $\CC$. Эта функция $v^{\bal}$ называется выметанием функции $u$ из открытого множества $\CC\setminus S$ на систему лучей $S$.
\end{theoB}

Очевидно, система из четырёх лучей $\RR^+$, $-\RR^+$, $i\RR^+$, $-i\RR^+$  $p$-допустима, если  и только если $p<2$. Отметим особенности этого частного случая.  
\begin{theoB}\label{BiR} Если  $M\in \sbh_*(\CC)$ кон\-е\-ч\-н\-о\-го типа при порядке $1$ с мерой Рисса $\upmu :=\frac1{2\pi}\Delta M\in \Meas^+(\CC)$, то найдутся   функция 
$M_{\RR}\in \sbh_*(\CC)$ конечного типа $\type_1[M_{\RR}]<+\infty$ при порядке $1$ с мерой Рисса 
\begin{equation}\label{nubalR}
\frac{1}{2\pi}\Delta M_{\RR}=:\upmu_{\RR}\in \Meas^+(\RR),
\quad \type_1[\upmu_{\RR}]<+\infty, 
\end{equation} 
удовлетворяющей условию Линделёфа рода $1$, а также 
функция $M_{i\RR}$ конечного типа 
$\type_1[M_{i\RR}]<+\infty$ при порядке $1$ с мерой Рисса \begin{equation}\label{nubaliR}
\frac{1}{2\pi}\Delta M_{i\RR}=:\upmu_{i\RR}\in \Meas^+(i\RR),
\quad \type_1[\upmu_{i\RR}]<+\infty, 
\end{equation} 
удовлетворяющей условию Линделёфа рода $1$, со свойствами  
\begin{subequations}\label{vbalM}
\begin{align}
M(iy)&=M_{i\RR}(iy)+M_{\RR}(iy)\quad\text{для каждого  $y\in \RR$,}
\tag{\ref{vbalM}i}\label{{vbalM}i}\\
M(x)&=M_{i\RR}(x)+M_{\RR}(x)\quad\text{для каждого $x\in \RR$,}
\tag{\ref{vbalM}r}\label{{vbalM}r}\\
M(z)&\leq M_{i\RR}(z)+M_{\RR}(z)\quad\text{для каждого $z\in \CC$.}
\tag{\ref{vbalM}z}\label{{vbalM}z}
\end{align}
\end{subequations}
При этом $\upmu\curlyeqprec_{|\Re|}\upmu_{\RR}$ в смысле определения\/ {\rm \ref{ab}}. 
\end{theoB}
\begin{proof} По теореме B\ref{BR} с системой четырёх лучей  $S:=\RR\cup i\RR$ в сочетании с теорема Вейерштрасса\,--\,Адамара\,--\,Линделёфа\,--\,Брело существует функция $M^{\bal}\in \sbh_*(\CC)$ конечного типа при порядке $1$ с мерой Рисса $\frac{1}{2\pi}\Delta M^{\bal}=:\upmu^{\bal}\in \Meas^+(\RR\cup i\RR)$ конечной верхней плотности при порядке $1$, удовлетворяющей условию Линделёфа рода $1$, и при этом 
\begin{subequations}\label{vbalMR}
\begin{align}
M^{\bal}(iy)&=M(iy)\quad\text{для всех $y\in \RR$,}
\tag{\ref{vbalMR}i}\label{{vbalMR}i}\\
M^{\bal}(x)&=M(x)\quad\text{для всех $x\in \RR$,}
\tag{\ref{vbalMR}r}\label{{vbalMR}r}\\
M^{\bal}(z)&\leq M(z)\quad\text{для всех $z\in \CC$.}
\tag{\ref{vbalMR}z}\label{{vbalMR}z}
\end{align}
\end{subequations}
Сужение $\upmu_{\RR}:=\upmu^{\bal}\bigm|_{\RR}\in \Meas^+(\RR)$ меры $\upmu$
на $\RR$, очевидно, также конечного порядка при порядке $1$ и удовлетворяет  условию Линделёфа \eqref{con:Lp+} рода $1$. 
Вновь по теореме Вейерштрасса\,--\,Адамара\,--\,Линделёфа\,--\,Брело существует функция $M_{\RR}\in \sbh_*(\CC)$ конечного типа $\type_1[M_{\RR}]<+\infty$ при порядке $1$ с мерой Рисса $\upmu_{\RR}$. Тогда функция $M_{i\RR}:=M^{\bal}-M_{\RR}\in \sbh_*(\CC)$ тоже конечного типа при порядке $1$, но с мерой Рисса 
$\upmu_{i\RR}=\upmu^{\bal}-\upmu_{\RR}\in \Meas^+(i\RR)$ конечной верхней плотности при порядке $1$, удовлетворяющей  условию Линделёфа рода $1$.  При этом 
соотношения \eqref{{vbalMR}i}, \eqref{{vbalMR}r}, \eqref{{vbalMR}z} влекут за собой соответственно соотношения \eqref{{vbalM}i}, \eqref{{vbalM}r}, \eqref{{vbalM}z}. Ввиду совпадения \eqref{{vbalM}i} функций $M$ и $M_{\RR}+M_{i\RR}$ на $i\RR$ по определению \eqref{fK:abp+} логарифмического интеграла для каждого интервала $(r,R]\subset \RR^+$ имеет  место равенство  $J_{i\RR}(r,R; M)=J_{i\RR}(r,R; M_{\RR}+M_{i\RR})$.
Отсюда по соотношению  \eqref{{Jll}m} предложения \ref{lemJl} получаем 
$\upmu \overset{\eqref{lJMmul}}{\curlyeqprec}_{|\Re|}\upmu_{\RR}+\upmu_{i\RR}$. Но $\supp \upmu_{i\RR}\subset i\RR $ и  по определению логарифмической субмеры интервалов \eqref{df:dDlLm}  и определению \ref{ab}  имеем 
$\upmu_{i\RR}\curlyeqprec_{|\Re|}0$, что влечёт за собой требуемое $\upmu \curlyeqprec_{|\Re|}\upmu_{\RR}$. 
\end{proof}

\section{Двустороннее выметание на мнимую ось}\label{bal_iR}
\setcounter{equation}{0}

Утверждения этого параграфа выводятся в основном из   результатов второй части \cite{KhII} нашей работы применительно  к выметанию рода $1$ на систему из двух лучей $\{i\RR^+, -i\RR^+ \}$. Такое выметание представляет собой один из ключевых этапов доказательства основной теоремы. 
При этом система двух лучей $\{i\RR^+, -i\RR^+ \}$ как точечное множество отождествляется с мнимой осью $i\RR$.  В частности, когда речь идет о выметании на пару лучей  $\{i\RR^+, -i\RR^+ \}$, то говорим о выметании на  $i\RR$. В целом нельзя утверждать, что процедура выметания рода $1$ на мнимую ось $i\RR$ --- это простое сочетание двух выметаний рода $1$:  отдельно из правой полуплоскости $\CC_{\rght}$  и из левой полуплоскости $\CC_{\lft}$. Так, 
для сохранения основных свойств заряда $\upnu$ при выметании отдельно из  
$\CC_{\rght}$  и из  $\CC_{\lft}$ требуется условие Бляшке 
\eqref{Blclq}--\eqref{Bl}  рода $1$ в  $\CC_{\rght}$  и в   $\CC_{\lft}$, а основные требуемые свойства выметания  рода $1$ для заряда $\upnu$ на мнимую ось $i\RR$ выполняются уже при двустороннем условии Бляшке из определения   \ref{df:bB}, которое по части [B\ref{B2}] предложения  \ref{pr2}  строго слабее. Тем не менее,  значительную часть результатов настоящего параграфа о двустороннем выметании удаётся свести к результатам о выметании рода $q=1$ из верхней полуплоскости $\CC^{\up}$ на $\RR$, изложенных в \cite[3.2, 4.1]{KhII}, а также в  \cite[\S~3]{kh91AA}, \cite[гл.~II]{KhDD92}.

\subsection{Двустороннее выметание рода $1$  заряда  на мнимую ось}

Для заряда $\upnu\in \Meas(\CC)$ наряду со считающими функциями 
\eqref{df:nup} используем как  {\it функцию распределения\/} $\upnu^{\RR}$ сужения $\upnu\bigm|_{\RR}$ заряда $\upnu$ {\it на вещественной оси\/}  $\RR$ \cite[(1.9)]{KhI}, определённую при $-\infty<x_1<x_2<+\infty$ равенствами 
\begin{equation}\label{nuR} 
\upnu^{\RR}(x_2)-\upnu^{\RR}(x_1):=\upnu\bigm|_{\RR}\bigl((x_1,x_2]\bigr), \quad
\upnu^{\RR}(0):=0,
\end{equation}
так и аналогичную её  {\it функцию распределения\/} $\upnu^{i\RR}$ сужения $\upnu\bigm|_{i\RR}$ заряда $\upnu$ {\it на мнимой оси\/} $i\RR$, определённую при $-\infty<y_1<y_2<+\infty$ по правилу
\begin{equation}\label{nuiR} 
\upnu^{i\RR}(y_2)-\upnu^{i\RR}(y_1):=\upnu\bigm|_{i\RR}\bigl((y_1,y_2]\bigr), \quad
\upnu^{i\RR}(0):=0.
\end{equation}
По построению \eqref{nuR} и \eqref{nuiR}  функции $\upnu^{\RR}\colon \RR\to \RR$ и $\upnu^{i\RR}\colon \RR \to \RR$ --- функции локально ограниченной вариации  на $\RR$, а в случае {\it меры\/} $\upnu \in \Meas^+(\CC)$ обе эти функции возрастающие на $\RR$.

Мы напоминаем и адаптируем основные понятия и утверждения из \cite{KhII},
 касающиеся выметания конечного рода $q$ из полуплоскости, применительно к правой и левой полуплоскостям в случае $q=1$.

Пусть $I:=(iy_1,iy_2]:=i(y_1,y_2]$ ---  интервал на $i\RR$, $-\infty<y_1<y_2<+\infty$. {\it Двустороннюю гармоническую меру для\/ $\CC\setminus i\RR$ в точке\/ $z\in \CC$} обозначаем как  функцию интервалов 
\begin{equation}\label{omega}
\omega \bigl(z,(iy_1,iy_2]\bigr)\overset{\text{\cite[3.1]{KhI}}}{:=}\omega_{\CC\setminus i\RR}(z,I):=\frac1{\pi}
\int_{y_1}^{y_2}\Bigl|\Re \frac{1}{iy-z}\Bigr| \dd y ,
\end{equation} 
равную делённому на $\pi$ углу, под которым виден интервал $I$ из точки $z\in \CC$ \cite[(3.1)]{kh91AA}, \cite[1.2.1, 3.1]{KhI}.
Если \textit{мера\/}  $\upmu \in \Meas^+(\CC)$ конечного типа удовлетворяет классическому условию Бляшке \eqref{Blcl}--\eqref{lchar} в $\CC_{\rght}$
и в $\CC_{\lft}$, то  определено классическое выметание $\upmu^{\bal}$ (рода $0$ \cite[определение 3.1]{KhII}) этой меры на мнимую ось \cite[следствие 4.1, теорема 4]{KhI} с носителем $\supp \upmu^{\bal} \subset i\RR $. Эту выметенную меру 
$\upmu^{\bal}$  можно задать через функцию распределения $(\upmu^{\bal})^{i\RR}$ из 
\eqref{nuiR}  меры $\upmu^{\bal}$   как 
\begin{equation}\label{mubal}
(\upmu^{\bal})^{i\RR}(y_2)
-(\upmu^{\bal})^{i\RR}(y_1)
\overset{\eqref{omega}}{:=}
\int\limits_{\CC} \omega\bigl(z, (iy_1,iy_2]\bigr) \dd \upmu(z)
\end{equation}
с нормировкой $\quad (\upmu^{\bal})^{i\RR}(0):=0$.

Гармонический заряд рода $1$ для $\CC_{\rght}$ в точке $z\in \CC_{\overline \rght}\setminus \{0\}$ 
определяется на интервалах $(iy_1,iy_2]\subset i\RR$, $-\infty<y_1<y_2<+\infty$, через \cite[определение 2.1]{KhII}    как функция интервалов
 \begin{multline}\label{Ocr}
\Omega_{\rght}\bigl(z,(iy_1,iy_2]\bigr):=\Omega_{\CC_{\rght}}^{[1]}\bigl(z,(iy_1,iy_2]\bigr)\\
:=\omega_{\rght}\bigl(z,(iy_1,iy_2]\bigr)-\frac{y_2-y_1}{\pi}\Re\frac{1}{z}, \quad z\neq 0,
\end{multline} 
где $\omega_{\rght}\bigl(z,(iy_1,iy_2]\bigr)$ --- классическая гармоническая мера 
интервала $(iy_1,iy_2]$ для  $\CC_{\rght}$ в точке $z$, равная 
делённому на $\pi$ углу, под которым виден интервал $(iy_1,iy_2]$ из точки 
$ z\in \CC_{\overline \rght}$  \cite[определение 4.3.1]{Rans}, \cite[\S~3]{KhI}.
Подобным же образом определяется гармонический заряд рода $1$ для $\CC_{\lft}$ в точке $z\in \CC_{\overline \lft}\setminus \{0\}$, а именно: 
$\Omega_{\lft}(z,\cdot):=\Omega_{\rght}(-\bar z,\cdot)$, 
$\omega_{\lft}(z,\cdot):=\omega_{\rght}(-\bar z,\cdot)$. 
 {\it Двусторонний гармонический заряд рода $1$ для\/ $\CC\setminus i\RR$ в точке\/ $z\in \CC$} определяем на интервалах $(iy_1,iy_2]$, $-\infty<y_1<y_2<+\infty$, как функцию интервалов  \cite[определение 3.1]{kh91AA}, \cite[определение 2.1.1]{KhDD92}, \cite[определение 2.1]{KhII} 
\begin{multline}\label{se:HC+1}
\Omega\bigl(z,(iy_1,iy_2]\bigr)
:=:\Omega_{\CC\setminus i\RR}^{[1]}\bigl(z,(iy_1,iy_2]\bigr):=
\begin{cases}
\Omega_{\rght}\bigl(z,(iy_1,iy_2]\bigr)&\text{ при $z\in \CC_{\overline \rght}$,}\\
\Omega_{\lft}\bigl(z,(iy_1,iy_2]\bigr)&\text{ при $z\in \CC_{\overline \lft}$}
\end{cases}
\\
\overset{\eqref{omega}}{=}\omega\bigl(z,(iy_1, iy_2]\,\bigr)
-\frac{y_2-y_1}{\pi}  \Bigl|\Re \frac{1}{z}\Bigr| \, , \quad z\neq 0,
\end{multline}
где двусторонняя гармоническая мера  
\begin{equation*}
\omega(z,\cdot):=\begin{cases}
\omega_{\rght}(z,\cdot)\quad&\text{при $z\in \CC_{\overline \rght}$,}\\
\omega_{\lft}(z,\cdot)\quad&\text{при $z\in \CC_{\overline \lft}$}
\end{cases}
\end{equation*}
для $\CC\setminus i\RR$ в точке $z\in \CC$ именно та, что определена в \eqref{omega}. 
 
Для  заряда $\upnu\in \Meas(\CC)$  конечного типа $\type[\upnu]<+\infty$ 
в  \cite[определение 3.1, теорема 1, замечание 3.3]{KhII} определялось 
выметание рода $1$ заряда $\upnu$ из правой полуплоскости $\CC_{\rght}$ на замкнутую левую полуплоскость $\CC_{\overline \lft}$.  Сужение 
такого  выметания  на мнимую ось $i\RR$ можно однозначным образом определить   через функцию распределения вида \eqref{nuiR} локально ограниченной вариации на $i\RR$ по правилу, согласованному с  \eqref{mubal}, как 
\begin{equation}\label{df:nurh}
\begin{split}
(\upnu^{\Bal}_{\rght})^{i\RR}(y_2)-(\upnu^{\Bal}_{\rght})^{i\RR}(y_1)
&=\int_{\CC_{\rght}\cap \DD} \omega_{\rght} \bigl(z, (iy_1,iy_2]\bigr)\dd \upnu (z)\\
+\int_{\CC_{\rght}\setminus \DD} \Omega_{\rght} \bigl(z, (iy_1,iy_2]\bigr)\dd \upnu (z)
&+\upnu\bigl((iy_1,iy_2]\bigr),\;  -\infty<y_1<y_2<+\infty,
\end{split}
\end{equation}
с возможной нормировкой $(\upnu^{\Bal}_{\rght})^{i\RR}(0):=0$. 
Аналогичным образом сужение на $i\RR$ выметания рода $1$ заряда $\upnu$ из левой полуплоскости $\CC_{\lft}$ на замкнутую правую полуплоскость $\CC_{\overline \rght}$ определяется   через функцию распределения  
\begin{equation}\label{df:blh}
(\upnu^{\Bal}_{\lft})^{i\RR}(y_2)-(\upnu^{\Bal}_{\lft})^{i\RR}(y_1)
:=\bigl((\upnu_{\leftrightarrow})^{\Bal}_{\rght}\bigr)^{i\RR}(y_2)-\bigl((\upnu_{\leftrightarrow})^{\Bal}_{\rght}\bigr)^{i\RR}(y_1),
\end{equation}
где заряд $\upnu_{\leftrightarrow}$ зеркально симметричный заряду $\upnu$
относительно $i\RR$ из \eqref{nus}.

Исходя из этих уже рассмотренных в \cite[\S~3]{KhI} видов выметания, можно   
определить, следуя  \cite[определение 3.2]{kh91AA}, \cite[определение 2.1.2]{KhDD92},  {\it двустороннее выметание\/ $\upnu^{\Bal}\in \Meas (i\RR)$ заряда $\upnu$ из $\CC\setminus i\RR$ на $i\RR$\/} как заряд $\upnu^{\Bal}\in \Meas (\CC)$ 
с носителем $\supp \upnu^{\Bal} \subset i\RR$ через функцию распределения  
\begin{equation}\label{df:2b}
\begin{split}
(\upnu^{\Bal})^{i\RR}(y_2)-(\upnu^{\Bal})^{i\RR}(y_1)&:=
(\upnu_{\rght}^{\Bal})^{i\RR}(y_2)-(\upnu_{\rght}^{\Bal})^{i\RR}(y_1)\\
+(\upnu^{\Bal}_{\lft})^{i\RR}(y_2)-(\upnu^{\Bal}_{\lft})^{i\RR}(y_1)&-
\upnu\bigl((iy_1,iy_2]\bigr),\quad  -\infty<y_1<y_2<+\infty.
\end{split}
\end{equation}
В более явном виде через двусторонний гармонический заряд \eqref{se:HC+1} рода $1$ двустороннее выметание $\upnu^{\Bal}\in \Meas (i\RR)$ заряда $\upnu$можно определить как 
 \begin{equation}\label{df:baliR}
\begin{split}
(\upnu^{\Bal})^{i\RR}(y_2)-(\upnu^{\Bal})^{i\RR}(y_1)&\overset{\eqref{df:2b}}{:=}
\int_{(\CC\setminus i\RR)\cap \DD} \omega\bigl(z, (iy_1,iy_2]\bigr)\dd \upnu(z)\\
+\int_{(\CC\setminus i\RR)\setminus \DD} \Omega \bigl(z, (iy_1,iy_2]\bigr)\dd \upnu(z)&+ \upnu \bigl((iy_1,iy_2]\bigr),
\end{split}
\end{equation} 
где $-\infty<y_1<y_2<+\infty$.
Если определить  \textit{$i\RR$-симметризацию заряда\/ $\upnu$,} исходя из определения \eqref{nus}, как заряд из \eqref{reprnu} вида 
\begin{equation}\label{nusym}
\upnu_{\rightleftarrows}\overset{\eqref{nus}}{:=}\frac{1}{2}(\upnu+\upnu_{\leftrightarrow})
\end{equation}
то по  \eqref{df:nurh}, \eqref{df:blh}, \eqref{df:2b}, \eqref{df:baliR} двустороннее выметание $\upnu^{\Bal}$ совпадает с сужением на мнимую ось выметания $i\RR$-симметризации \eqref{nusym} заряда $\upnu$:
\begin{equation}\label{nust}
\upnu^{\Bal}=(\upnu_{\rightleftarrows})_{\rght}^{\Bal}\bigm|_{i\RR}.
\end{equation}

\begin{theoB}\label{theoB1} Пусть $\upnu\in \Meas(\CC)$  --- заряд конечной верхней плотности ${\type_1}[\upnu] <+\infty$.  Тогда существует двустороннее выметание рода $1$ на мнимую ось $\upnu^{\Bal}\in \Meas(i\RR)$  из  \eqref{df:2b}--\eqref{df:baliR}, для которого 
\begin{subequations}\label{iRbnu}
\begin{align}
\bigl|\upnu^{\Bal}\bigr|^{\rad}(r)&=O(r\ln r) \quad\text{при  $r\to +\infty$},
\tag{\ref{iRbnu}l}\label{iRbnui}
\\
\intertext{а  при $0\notin \supp \upnu$ 
имеет место соотношение}
\bigl|\upnu^{\Bal}\bigr|^{\rad}(r)&=O(r^2)\quad\text{при  $0<r\to 0$}.
\tag{\ref{iRbnu}o}\label{iRbnu0}
\end{align}
\end{subequations} 

Пусть   заряд $\upnu\in \Meas (\CC)$ с $\type_1[\upnu]<+\infty$ 
удовлетворяет ещё и двустороннему условию Бляшке \eqref{cB2} рода $1$ вне $i\RR$. Тогда
\begin{enumerate}
\item[{\rm [{F}]}] 
 $\upnu^{\Bal}$ --- заряд конечной верхней плотности
${\type_1}[\upnu^{\Bal}] <+\infty$. 

\item[{\rm [{S}]}]  
Если  заряд $\upnu$ отделён углами от $i\RR$,
то в обозначении \eqref{df:nup} имеем
\begin{equation}\label{trnuair}
\sup_{y\in \RR}\sup_{t\in (0,1]}\frac{|\upnu^{\Bal}\bigr|(iy, t)}{t}<+\infty.
\end{equation}

\item[{\rm [{L}]}]  Разность зарядов  $\upnu-\upnu^{\Bal}$ 
удовлетворяет условию Линделёфа \eqref{con:Lp+} рода $1$. Так, 
если для заряда $\upnu$ выполнено условие \eqref{con:LpIm} или
условие Линделёфа \eqref{con:Lp+} рода $1$, то и выметание $\upnu^{\Bal}\in \Meas(i\RR)$ удовлетворяют условию Линделёфа рода $1$.
\end{enumerate}
\end{theoB}
\begin{proof} Существование двустороннего выметания $\upnu^{\Bal}$ и соотношение  \eqref{iRbnui} согласно представлению \eqref{nust} --- частный случай  \cite[теорема 3, п.~2]{KhII} и установлено ещё в \cite[лемма 2.1.2, (1.16)]{KhDD92}. Из того же представления  \eqref{nust} соотношение \eqref{iRbnu0} --- частный случай 
\cite[теорема 1]{KhII} и имеется в  \cite[лемма 2.1.2, (1.17)]{KhDD92}. Часть [F] 
согласно  представлению \eqref{nust} --- частный случай  \cite[теорема 3, п.~4]{KhII}
и получена ещё в \cite[теорема 3.1]{kh91AA}, \cite[теорема 2.2.1]{KhDD92}.
Часть [S] согласно  представлению \eqref{nust} --- частный случай 
\cite[следствие 3.1, п.~(ii), (3.24)]{KhII}, который в неявной форме имеется в  \cite[теорема 3.3]{kh91AA}, а в явном виде  --- в  \cite[теорема 2.2.2]{KhDD92}. 
Заключительная часть [L] при дополнительном условии   \eqref{con:LpIm} была доказана  в \cite[теорема 3.2]{kh91AA}. Это доказательство почти дословно можно адаптировать под доказательство более общего утверждения  из [L] о выполнении условия Линделёфа для разности зарядов $\upnu-\upnu^{\Bal}$ без условия   \eqref{con:LpIm} или
условия Линделёфа \eqref{con:Lp+} рода $1$ для $\upnu$. 
Последнее без всяких дополнительных условий на заряд $\upnu$ конечной верхней плотности при порядке $1$  детально проделано в диссертации  \cite[теорема 2.3.1]{KhDD92}.
\end{proof}

\subsection{Двустороннее выметание рода\/ $1$ $\updelta$-субгармонической функции на мнимую ось}

Пусть $v\in \dsbh_*(\CC)$. Функцию $v^{{\Bal}}\in \dsbh_*(\CC)$ называем {\it двусторонним  выметанием функции\/ $v$ из\/ $\CC\setminus i\RR$ на $i\RR$}, если $v^{{\Bal}}=v$ на\/ $i\RR$ вне полярного множества  и сужение  $v^{{\Bal}}\bigm|_{\CC\setminus i\RR}$ --- гармоническая функция на   $\CC\setminus i\RR$ \cite[определение 4.1]{KhII}. 
Сочетание теоремы  Вейерштрасса\,--\,Адамара\,--\,Линделёфа\,--\,Брело 
с теоремой B\ref{theoB1} вместе с \cite[теоремы 6,7]{KhII} и \cite[теорема 2.1.1]{KhDD92} позволяет дать следующую сводку. 

\begin{theoB}
\label{Balv} Пусть $v\in \dsbh_*(\CC)$ --- функция c зарядом Рисса $\upnu$ ко\-н\-е\-ч\-н\-ой верхней плотности 
$\type_1[\upnu ]<+\infty$  при порядке $1$ и функция $v$ представима в виде разности субгармонических функций
\begin{equation}\label{reprv}
v:=v_+-v_-, \quad v_{\pm}\in \sbh_*(\CC), \quad \ord[v_{\pm}]\overset{\eqref{ord}}{\leq} 1; 
\quad\type_1[\upnu]<+\infty. 
\end{equation} 
 Тогда существует двустороннее выметание\/ $v^{\Bal}\in \dsbh_*(\CC)$ из\/ $\CC\setminus i\RR$ на $i\RR$ c зарядом Рисса $\upnu^{\Bal}\in \Meas(i\RR)$ с функцией распределения \eqref{df:2b}--\eqref{df:baliR}  представимое в виде разности субгармонических функций
\begin{equation}\label{reprvB}
v^{\Bal}:=u_+-u_-, \quad u_{\pm}\in \sbh_*(\CC), \quad \ord[u_{\pm}]\overset{\eqref{ord}}{\leq}  1; \quad \frac{1}{2\pi}\Delta v:=\upnu^{\Bal}.
\end{equation} 
Пусть в  \eqref{reprv} функции $v_{\pm}$ конечного типа ${\type_1}[v_{\pm}]<+\infty$. Тогда 
\begin{enumerate}
\item[{\rm [{F}]}] если  заряд  $\upnu$ удовлетворяет двустороннему условию Бляшке \eqref{cB2} рода $1$, то две функции $u_{\pm}$ в представлении \eqref{reprvB} можно выбрать так, что $\type_1[{\mathsf  C}_{u_{\pm}}]\overset{\eqref{typevf}}{<}+\infty$;  

\item[{\rm [{L}]}] если заряд  $\upnu$ удовлетворяет двустороннему условию Бляшке \eqref{cB2} рода $1$ и  условию Линделёфа рода $1$ из \eqref{con:Lp+}, то две функции $u_{\pm}$ в \eqref{reprvB} можно выбрать конечного типа  $\type_1[u_{\pm}]<+\infty$;
\item[{\rm [{S}]}]
Если  при этом заряд $\upnu$ еще и {отделён углами от мнимой оси}  $i\RR$, то для некоторых $C\in \RR^+$ и $y_0\in \RR^+$ имеем неравенства 
\begin{equation*}
u_\pm (iy)\overset{\eqref{df:MCBc}}{\leq}\mathsf{C}_{u_\pm} (iy,1) \overset{\eqref{trnuair}}{\leq} u_\pm (iy)+ C\quad \text{при всех $|y|\geq y_0$}.
\end{equation*}
\end{enumerate}
\end{theoB}

\section{Две   конструкции с  логарифмическими функциями интервалов\\ и выметаниями мер на мнимую ось}
\setcounter{equation}{0}

\begin{propos}[{для последовательностей см. \cite[лемма 3.1]{MR}, \cite[лемма 22.2]{RC}, \cite[лемма 1.1.]{Kha89},  \cite[лемма 1]{kh91AA}, для зарядов  частично  --- в \cite[лемма 2.4.2]{KhDD92}}]\label{addsec}
Пусть для заряда $\eta \in \Meas(\CC)$ имеем  $0\not\in \supp \eta$ и 
\begin{equation}\label{lnu}
\begin{split}
M^{\rght}:=\sup_{0< r<R<+\infty} l_{\eta}^{\rght}(r,R)&< +\infty\\
\Bigl(\text{соответственно } M^{\lft}:=\sup_{0< r<R<+\infty}l_{\eta}^{\lft}(r,R)&< +\infty \Bigr).
\end{split}
\end{equation}
Тогда определена возрастающая функция 
\begin{equation}\label{ainf}
\begin{split}
a(t)&\overset{\eqref{df:dDlL1}}{:=}-\sup_{s\geq t} l_{\eta}^{\rght}(s)
=\inf_{s\geq t} \bigl(-l_{\eta}^{\rght}(s)\bigr),
\quad t\in \RR^+ \\
\Bigl(\text{соответственно } a(t)&\overset{\eqref{df:dDlL1}}{:=}-\sup_{s\geq t} l_{\eta}^{\lft}(s)
=\inf_{s\geq t} \bigl(-l_{\eta}^{\lft}(s)\bigr), \quad t\in \RR^+ \Bigr),
\end{split}
\end{equation}
однозначно определяющая  меру  $\alpha \in \Meas^+(\RR^+)$ 
с\/ $0\neq \supp \alpha$  
через её  функцию распределения 
\begin{subequations}\label{adist}
\begin{align}
\alpha^{\RR}(x)&\overset{\eqref{nuR}}{:=}\int_0^x t \dd a(t)
=xa(x)-\int_0^xa(t)\dd t, \quad x\in \RR^+
\tag{\ref{adist}r}\label{{adist}r}
 \\
\intertext{{\rm \Large(}соответственно} \alpha^{\RR}(x)&\overset{\eqref{nuR}}{:=}-\int_0^{-x} t \dd a(t)
=xa(-x)+\int_0^{-x}a(t)\dd t, \quad x\in -\RR^+\Bigr),
\tag{\ref{adist}l}\label{{adist}l}
\end{align}
\end{subequations}
и с этой мерой $\alpha$ выполнено соотношение  
\begin{equation}\label{|l|}
\begin{split}
\sup_{0< r<R<+\infty} &\bigl|l_{\eta+\alpha}^{\rght}(r,R)\bigr|\leq 2M^{\rght}< +\infty\\
\Bigl(\text{соответственно }\sup_{0< r<R<+\infty} &\bigl|l_{\eta+\alpha}^{\lft}(r,R)\bigr|
\leq 2M^{\lft}< +\infty\Bigr).
\end{split}
\end{equation}
Если заряд $\eta$ конечной верхней плотности $\type_1[\eta]<+\infty$, 
то мера $\alpha$ из   \eqref{ainf}--\eqref{adist} также конечной верхней плотности  $\type_1[\alpha]<+\infty$.
\end{propos}
\begin{proof} Достаточно рассмотреть случай с $l_{\eta}^{\rght}$.  Возрастание фу\-н\-к\-ции $a$ очевидно по её построению \eqref{ainf}. Из её возрастания следует, что функция распределения   $\alpha^{\RR}$ из \eqref{{adist}r} также возрастающая. Следовательно, заряд $\alpha\in \Meas(\RR^+)$ с такой функцией распределения
$\alpha^{\RR}$  --- это положительная мера на $\RR^+$. Ввиду $0\notin \supp \eta$ из построения \eqref{ainf} функция $a$ постоянна в некоторой окрестности нуля, функция распределения $\alpha^{\RR}$ --- тождественный  нуль  в этой же окрестности  и $0\notin \supp \alpha$.  Кроме того, по построению 
\eqref{ainf} для всех $t\in \RR^+$ имеем 
\begin{equation}\label{aup}
a(t)\overset{\eqref{ainf}}{=}\inf_{s\geq t}\bigl( -l_{\eta}^{\rght}(s)\bigr)\leq 
 -l_{\eta}^{\rght}(t), \quad a(0)\leq 0, 
\end{equation}
а также при всех $t\in \RR^+$
\begin{equation}\label{a}
 a(t)\overset{\eqref{ainf}}{=}  -\sup_{s\geq t} l_{\eta}^{\rght}(s)
\geq  -\sup_{0< r<R<+\infty} l_{\eta}^{\rght}(r,R)\overset{\eqref{lnu}}{=} -M^{\rght}.
\end{equation}
Кроме того,  по определению \eqref{ainf} при всех $t\in \RR^+$ имеем 
\begin{equation*}
l_{\alpha}^{\rght} (t)\overset{\eqref{adist}}{=}
\int_0^t\frac{1}{x}\dd \alpha^{\RR}(x)\overset{\eqref{adist}}{:=}\int_0^t  \dd a(t)
=a(t)-a(0),
\end{equation*}
откуда
\begin{equation}\label{nua}
l_{\eta+\alpha}^{\rght}(t)
=l_{\eta}^{\rght}(t)+a(t)-a(0)\overset{\eqref{aup}}{\leq}-a(0)
\overset{\eqref{a}}{\leq}M^{\rght}, 
\end{equation}
а также
\begin{multline*}
l_{\eta+\alpha}^{\rght}(t) \overset{\eqref{nua}}{=}
l_{\eta}^{\rght}(t)+a(t)-a(0)
\overset{\eqref{ainf}}{=}
\inf_{s\geq t}\bigl(l_{\eta}^{\rght}(t) -l_{\eta}^{\rght}(s)\bigr)-a(0)\\
=-\sup_{s\geq t}\bigl(l_{\eta}^{\rght}(s)-l_{\eta}^{\rght}(t) \bigl)-a(0)
=-\sup_{s\geq t}\bigl(l_{\eta}^{\rght}(s,t) \bigl)-a(0)\\
\overset{\eqref{aup}}{\geq}    -\sup_{0< r<R<+\infty} l_{\eta}^{\rght}(r,R)
\overset{\eqref{lnu}}{=} -M^{\rght}.
\end{multline*}
Отсюда и из \eqref{nua} получаем
\begin{equation*}
\sup_{t\in \RR^+} \bigl|l_{\eta+\alpha}^{\rght}(t)\bigr|\leq M^{\rght},
\end{equation*}
что влечёт за собой \eqref{|l|}:
\begin{equation*}
\sup_{0<r<R<+\infty} \bigl|l_{\eta+\alpha}^{\rght}(r,R)\bigr|\leq 
\sup_{R\in \RR^+} \bigl|l_{\eta+\alpha}^{\rght}(R)\bigr|+\sup_{r\in \RR^+} \bigl|l_{\eta+\alpha}^{\rght}(r)\bigr|\overset{\eqref{|l|}}{\leq} 
M^{\rght}.
\end{equation*}
Если $\eta$ --- заряд конечной верхней плотности, то 
\begin{equation}\label{l2}
\bigl|l_{\eta}^{\rght}(r,2r)\bigr|\leq \int_r^{2r}\frac{1}{t}\dd |\eta|^{\rad}(t)\leq
\frac{1}{r}|\eta|^{\rad}(2r)=O(1) \quad\text{при $r\to +\infty$}. 
\end{equation}
Отсюда и из \eqref{|l|}
\begin{multline*}
\frac{1}{2r}\bigl(\alpha^{\rad}(2r)-\alpha^{\rad}(r)\bigr)
\leq \int_r^{2r} \frac{1}{t}\dd \alpha^{\rad}(t)
=l_{\alpha}^{\rght}(r,2r)\\
\overset{\eqref{|l|}}{\leq} 
2M^{\rght}+\bigl|l_{\eta}^{\rght}(r,2r)\bigr|
\overset{\eqref{l2}}{=}O(1) \quad\text{при $r\to +\infty$.}
\end{multline*}
Это  означает, что мера $\alpha$ конечной верхней плотности
$\type_1[\alpha]<+\infty$.
\end{proof}

\begin{propos}\label{pr52} Пусть две меры $\upnu\in \Meas^+(\CC)$
и $\upmu\in \Meas^+(\CC)$ конечного типа 
$\type_1[\upnu+\upmu]<+\infty$ удовлетворяют условию \eqref{lJMmul}.
Тогда найдётся  мера $\alpha\in \Meas^+(\RR)$ конечной верхней плотности  $\type_1[\alpha]<+\infty$, для которой существует двустороннее выметание 
\begin{equation}\label{varheta}
(\upnu+\alpha-\upmu)^{\Bal}=:\vartheta=\vartheta^+-\vartheta^-\in \Meas(i\RR)
\end{equation} 
заряда $\upnu+\alpha-\upmu\in \Meas(\CC)$ из $\CC\setminus i\RR$ на $i\RR$ конечной верхней плотности $\type_1\bigl[|\vartheta|\bigr]<+\infty$, где 
$\vartheta^+$ и $\vartheta^-$ --- соответственно верхняя и нижняя вариации заряда $\vartheta$.
Если в обозначении   из  \eqref{nu-} для  меры $\vartheta^+_{\times}$, центрально симметричной мере $\vartheta^+$,  положим 
\begin{subequations}\label{bg}
\begin{align}
\beta&:=\vartheta^-+\vartheta^+_{\times}\in \Meas^+(i\RR), \quad\text{где, очевидно, } \type_1[\beta]<+\infty, 
\tag{\ref{bg}$\beta$}\label{bgb}
\\
\gamma&:=\vartheta^++\vartheta^+_{\times}\in \Meas^+(i\RR),\quad\text{где, очевидно, } \type_1[\gamma]<+\infty,
\tag{\ref{bg}$\gamma$}\label{bgg}
\end{align}
\end{subequations}
то мера $\gamma$ чётная,  имеет место равенство   
\begin{equation}\label{asg}
(\upnu+\alpha-\upmu)^{\Bal}+\beta=\vartheta+\beta =\gamma,
\end{equation}
мера $\gamma$ удовлетворяет условию Линделёфа,  а разность зарядов 
\begin{equation}\label{Lnumu}
(\upnu+\alpha-\upmu)-(\upnu+\alpha-\upmu)^{\Bal}
\end{equation}
удовлетворяет условию Линделёфа \eqref{con:LpIm} рода $1$. 
В частности,  и  сумма меры $\gamma$ из \eqref{asg} с зарядом \eqref {Lnumu}, равная 
\begin{equation}\label{sum}
\Bigl((\upnu+\alpha-\upmu)^{\Bal}+\beta\Bigr)
+\Bigl((\upnu+\alpha-\upmu)-(\upnu+\alpha-\upmu)^{\Bal}\Bigr)
=\upnu+\alpha+\beta-\upmu, 
\end{equation}
удовлетворяет условию Линделёфа. При этом
\begin{enumerate}
\item[{\rm [{L}]}]  если мера $\upmu$ удовлетворяет условию Линделёфа \eqref{con:LpIm} рода $1$, то мера $\upnu+\alpha+\beta$  также 
удовлетворяет условию Линделёфа  рода $1$.
\item[{\rm [{S}]}] если носители мер $\upnu$ и $\upmu$ не пересекаются с парой замкнутых вертикальных углов \eqref{anglei}, то в \eqref{bgg}--\eqref{asg} можно выбрать $\gamma:=c\lambda_{i\RR}$, где $c\in \RR^+$, а меру $\beta$ так, что выполнено  соотношение 
\begin{equation}\label{trnuairesg}
\sup_{y\in \RR}\sup_{t\in (0,1)}\frac{\beta(iy, t)}{t} <+\infty.
\end{equation}
\end{enumerate} 
 \end{propos}
\begin{proof} Не умаляя общности, можно считать, что 
\begin{equation*}
0\not\in \supp \upnu\cup \supp \upmu. 
\end{equation*}
Для заряда $\eta:=\upnu-\upmu\in \Meas(\CC)$ по предложению  \ref{addsec} существует мера $\alpha\in \Meas(\RR)$ конечной верхней плотности $\type_1[\alpha]<+\infty$, для которой  $\type_1[\eta+\alpha]<+\infty$ и выполнены оба соотношения из \eqref{|l|}. Это означает, что для заряда   $\eta+\alpha\in \Meas(\CC)$ выполнено   условие Бляшке \eqref{Bl} рода $1$ в обеих полуплоскостях $\CC_{\rght}$ и $\CC_{\lft}$. По  части [B\ref{B2}] предложения \ref{pr2} заряд $\eta+\alpha$ удовлетворяет двустороннему условию Бляшке     рода $1$ вне $i\RR$ из определения \ref{df:bB}. Следовательно, по части [F] теоремы B\ref{theoB1} существует  двустороннее выметание $(\eta+\alpha)^{\Bal}\in \Meas(i\RR)$ рода $1$ заряда $\eta+\alpha$ на мнимую ось конечной верхней плотности $\type_1\bigl[(\eta+\alpha)^{\Bal}\bigr]<+\infty$.

Представим заряд
\begin{equation}\label{varth}
\vartheta:=(\upnu +\alpha -\upmu)^{\Bal}=(\eta+\alpha)^{\Bal}\in \Meas^+(i\RR), \quad
\type_1[\vartheta]<+\infty, 
\end{equation}  в виде разности его положительной и отрицательной вариаций 
\begin{equation*}
\vartheta :=\vartheta^+-\vartheta^-, \quad 
\type\bigl[\vartheta^{\pm}\bigr]<+\infty,\quad \vartheta^{\pm}\in \Meas(i\RR).
\end{equation*}
Тогда 
\begin{equation}\label{vt}
(\upnu +\alpha -\upmu)^{\Bal}+\vartheta^-+\vartheta^+_{\times}=\vartheta^++\vartheta^+_{\times},
\end{equation}
где мера $\vartheta^+_{\times}\in \Meas^+(i\RR)$ центрально симметрична мере  $\vartheta^+\in \Meas^+(i\RR)$ в смысле \eqref{nu-}, а значит мера 
$\gamma :=\vartheta^++\vartheta^+_{\times}\in \Meas^+(i\RR)$ чётная и конечной верхней плотности $\type_1[\gamma]<+\infty$. Полагая 
$\beta:=\vartheta^-+\vartheta^+_{\times}\in \Meas^+(i\RR)$, 
из \eqref{vt}  получаем требуемое представление \eqref{asg}. При этом по части [L] теоремы B\ref{theoB1} разность зарядов \eqref{Lnumu} удовлетворяет условию Линделёфа.  
Очевидно, заряд из  \eqref{asg} удовлетворяет условию Линделёфа как чётный заряд. Сумма двух  зарядов  \eqref{asg} и \eqref{Lnumu}, удовлетворяющих условию Линделёфа \eqref{con:Lp+}, даёт заряд $\upnu+\alpha+\beta-\upmu$ из \eqref{sum}, также удовлетворяющий условию Линделёфа рода $1$.
Отсюда, если  мера $\upmu$  удовлетворяет условию Линделёфа  \eqref{con:Lp+}, то его сумма с зарядом $\upnu+\alpha+\beta-\upmu$ удов\-л\-е\-т\-воряет условию Линделёфа, что доказывает и часть [L] предложения \ref{pr52}.  

Если в условиях части [S] доказываемого предложения \ref{pr52} носители мер $\upnu$ и $\upmu$ не пересекаются с парой углов \eqref{anglei}, то по части [S] теоремы B\ref{theoB1} полная вариация $|\vartheta|$  заряда $\vartheta$ из \eqref{varth} удовлетворяет  условию вида   \eqref{trnuair}, которое можно продолжить на всю мнимую ось, а именно:    
\begin{equation*}
\sup_{y\in \RR}\sup_{t\in (0,1)}\frac{|\vartheta|(iy, t)}{t}<+\infty, 
\end{equation*}
что по построению для меры $\beta:=\vartheta^-+\vartheta^+_{\times}$ даёт требуемое в [S] второе соотношение из \eqref{trnuairesg}, а для меры $\gamma :=\vartheta^++\vartheta^+_{\times}$ имеем 
\begin{equation*}\label{trnuaivt}
c:=\sup_{y\in \RR}\sup_{t\in (0,1)}\frac{\gamma (iy, t)}{t}<+\infty. 
\end{equation*}
Из последнего сразу следует положительность меры $c\lambda_{i\RR}-\gamma\in \Meas^+(i\RR)$. Переопределим теперь меру  $\beta$ как сумму мер $\beta+(c\lambda_{i\RR}-\gamma)$. Тогда с мерой $c\lambda_{i\RR}$ в качестве меры $\gamma$  по-прежнему имеем равенство \eqref{asg} со всеми требуемыми выше в предложении \ref{pr52} свойствами для мер $\alpha \in \Meas(\RR)$, 
$\beta\in \Meas(i\RR)$, включая \eqref{trnuairesg}, и $\gamma :=c\lambda_{i\RR}\in \Meas^+(i\RR)$.   
\end{proof}
\section{Доказательство основной теоремы}\label{proof:2}
\setcounter{equation}{0} 

\subsection{Случай функции $M$, гармонической на   паре вертикальных углов \eqref{anglei}}
По условиям меры $\upnu\in \Meas^+(\CC)$ и $\upmu\in \Meas^+(\CC)$  конечной верхней плотности $\type_1[\upnu+\upmu]<+\infty$,  мера $\upnu$ отделена  углами  от мнимой оси и выполнено условие  \eqref{lJMmul}. По предположению  
функция $M\in \sbh_*(\CC)$ конечного типа $\type_1[M]<+\infty$ гармоническая  на  паре  углов \eqref{anglei}  с мерой Рисса $\frac{1}{2\pi}\Delta M\geq \upmu$. Таким образом, по предположению мера $\upmu$ также {\it отделена углами от мнимой оси.\/} Более того,  условие \eqref{lJMmul} не нарушится, если в нём заменить меру $\upmu$ на меру Рисса $\frac{1}{2\pi}\Delta M$ функции $M$. Следовательно, далее при доказательстве, не умаляя общности, можно считать, что  $\upmu=\frac{1}{2\pi}\Delta M$ --- мера Рисса функции $M$. Тогда, в частности, мера $\upmu$ удовлетворяет условию Линделёфа \eqref{con:Lp+} рода $1$.  
Таким образом,    выполнены посылки частей  [L] и [S]  предложения \ref{pr52}. Следовательно, по заключениям этих частей [L] и [S] существуют меры $\alpha\in \Meas^+(\RR)$ и $\beta, \gamma\in \Meas^+(i\RR)$ со всеми  указанными в предложении  \ref{pr52} свойствами \eqref{bg}--\eqref{sum}.  По заключению    [S] предложения \ref{pr52} вместе с \eqref{trnuairesg} их можно выбрать так, что 
\begin{equation}\label{asgl}
(\upnu+\alpha-\upmu)^{\Bal}+\beta\overset{\eqref{asg}}{=}
\gamma:=c\lambda_{i\RR}\quad\text{для некоторого числа  $c\in \RR^+$}.
\end{equation}
Пусть $u\in \sbh_*(\CC)$ с мерой Рисса $\frac{1}{2\pi}\Delta u=\upnu$ конечной верхней плотности $\type_1[\upnu]<+\infty$ из условия основной теоремы. По теореме Вейерштрасса\,--\,Адамара\,--\,Линделёфа\,--\,Брело 
через представление  \eqref{repru}
строятся функция $u_1\in \sbh_*(\CC)$ порядка 
$\ord[u_1]\overset{\eqref{typev}}{\leq} 1$ с мерой Рисса $\frac{1}{2\pi}\Delta u_1=\upnu$ и функция $v_1\in \sbh_*(\CC)$ порядка $\ord[v_1]\overset{\eqref{typev}}{\leq} 1$ с мерой Рисса $\alpha+\beta$, для которых 
$u_1+v_1\overset{\eqref{repru}}{=}u+H+v_1\in \sbh_*(\CC)$, где $H\in \har(\CC)$. 
По построению  функция $u+(v_1+H)=u_1+v_1$ порядка 
 $\ord[u+(v_1+H)]\leq 1$ с мерой Рисса $\upnu+\alpha+\beta$, 
которая по части [L] предложения \ref{pr52} удовлетворяет условию Линделёфа рода $1$. Следовательно, по  теореме Вейерштрасса\,--\,Адамара\,--\,Линделёфа\,--\,Брело  
функция $u+(v_1+H)\in \sbh_*(\CC)$ конечного типа $\type_1[u+(v_1+H)]<+\infty$.  

Рассмотрим $\updelta$-субгармоническую функцию $\bigl(u+(v_1+H)\bigr)-M\in \dsbh_*(\CC)$ конечного типа с зарядом Рисса $\upnu+\alpha+\beta-\upmu\in \Meas(\CC)$ конечной верхней плотности  при порядке $1$. По теореме  B\ref{Balv} для неё можно построить двустороннее выметание  \eqref{reprvB} 
$$
\Bigl(\bigl(u+(v_1+H)\bigr)-M\Bigr)^{\Bal}\in \dsbh_*(\CC)
$$ из $\CC\setminus i\RR$ на $i\RR$ так, что
\begin{equation}\label{uv1l}
\Bigl(\bigl(u+(v_1+H)\bigr)-M\Bigr)(iy)
=\Bigl(\bigl(u+(v_1+H)\bigr)-M\Bigr)^{\Bal}(iy)\quad\text{для  $y\in \RR$,}
\end{equation}
лежащих \textit{вне некоторого полярного множества из\/} $\RR$,
 с зарядом Рисса
\begin{equation}\label{beta}
(\upnu+\alpha+\beta-\upmu)^{\Bal}=(\upnu+\alpha-\upmu)^{\Bal}
+\beta\overset{\eqref{asg}}{=}\gamma\overset{{\rm [L]}}{=}
c\lambda_{i\RR}, 
\end{equation} 
где первое равенство следует из того, что   $\supp \beta\subset i\RR$.  
Из  равенств \eqref{beta} сразу вытекает, что функцию  $$\Bigl(\bigl(u+(v_1+H)\bigr)-M\Bigr)^{\Bal}\in \dsbh_*(\CC)$$
 можно отождествить на $\CC$ вне полярного множества с 
субгармонической функцией порядка не выше $1$ с мерой Рисса $c\lambda_{i\RR}$, которая   по  теореме Вейерштрасса\,--\,Адамара\,--\,Лин\-д\-е\-л\-ё\-фа\,--\,Брело   представляется для некоторых чисел $a,b\in \CC$ в виде явно выписываемой  субгармонической функции  
\begin{equation}\label{piv}
 \pi c|\Re z|+\Re(az+b)=\Bigl(\bigl(u+(v_1+H)\bigr)-M\Bigr)^{\Bal}(z), \quad z\in \CC.
\end{equation}  
конечного типа при порядке $1$. Если положить 
\begin{equation}\label{con:v}
v(z)=(v_1+H)(z)-\Re(az+b), \quad z\in \CC,
\end{equation}
то функция $u+v\in \sbh_*(\CC)$ конечного типа $\type_1[u+v]<+\infty$ при порядке $1$, а 
согласно равенствам   \eqref{piv} и  \eqref{uv1l} для всех $y\in \RR$, лежащих \textit{вне некоторого полярного множества,} имеем тождество 
\begin{equation}\label{uvm}
u(iy)+v(iy)-M(iy)\equiv \pi \Re (iy)\equiv 0 \quad \Longrightarrow
\quad u(iy)+v(iy)\equiv M(iy).
\end{equation} 
Отсюда в силу непрерывности функции $M$ на $i\RR$ и полунепрерывности сверху функции $u+v$ на $i\RR$ последнее тождество выполнено для всех $y\in \RR$. Таким образом,   для субгармонической функции $U:=u+v\in \sbh_*(\CC)$ конечного типа  $\type_1[U]=\type_1[u+v]<+\infty$ 
с мерой Рисса 
$$
\frac{1}{2\pi}\Delta U= \frac{1}{2\pi}\Delta u+\frac{1}{2\pi}\Delta v\geq \frac{1}{2\pi}\Delta u= \upnu
$$
 из  тождества \eqref{uvm} для всех $y\in \RR$
получаем \eqref{u=M}, что завершает доказательство части [sbh] в условиях 
гармоничности $M$ на паре углов \eqref{anglei}. 

Перейдём к доказательству части [Hol]. По построению \eqref{con:v} функция $v\in \sbh_*(\CC)$ конечного порядка $\ord[v]\leq 1$, а для части $\beta$ её   меры Рисса 
\begin{equation}\label{albv}
\frac{1}{2\pi}\Delta v\overset{ \eqref{con:v}}{=}\frac{1}{2\pi}\Delta v_1=\alpha +\beta,
\quad \alpha\in \Meas^+(\RR), \quad \beta\in \Meas^+(i\RR),
\end{equation}
 выполнено  условие \eqref{trnuairesg}. Функцию $v$ можно представить в виде суммы $v=v_{\alpha}+v_{\beta}$ двух функций порядка 
$\ord[v_{\alpha}]\leq 1$ и $\ord[v_{\beta}]\leq 1$ с мерами Рисса соответственно $\alpha\overset{\eqref{albv}}{\in} \Meas^+(\RR)$ и $\beta\overset{\eqref{albv}}{\in} \Meas^+(i\RR)$. 
Функции распределения $n_{\sf A}^{\RR}:=\lfloor\alpha^{\RR}\rfloor$  
и $n_{\sf B}^{i\RR}:=\lfloor\beta^{i\RR}\rfloor$ на $\RR$, заданные через целую часть $\lfloor\cdot \rfloor$
однозначно определяют соответственно последовательности точек ${\sf A}$ на вещественной оси и ${\sf B}$ на мнимой оси со считающими мерами 
$n_{\sf A}$ и  $n_{\sf B}$ конечного типа при порядке $1$. С помощью стандартной техники устанавливается
  \begin{lemma}[{см. \cite[лемма 2.2]{kh91AA}, \cite[лемма 2.4.3]{KhDD92}}] Существует целая функция $f_{\sf A}$ с последовательностью нулей $\Zero_{f_{\sf A}}={\sf A}\subset \RR$, для которой 
\begin{equation}\label{fA}
\sup_{|y|\geq 1} \Bigl|\ln \bigl|f_{\sf A}(iy)\bigr|-v_{\alpha}(iy)\Bigr|<+\infty, 
\end{equation}
а при условии \eqref{trnuairesg} на меру $\beta$ существует функция $f_{\sf B}$ с последовательностью нулей $\Zero_{f_{\sf B}}={\sf B}\subset i\RR$, для которой 
\begin{equation}\label{fB}
\ln |f_{\sf B}|\leq v_{\beta}\quad \text{при всех $y\in \RR$.}
\end{equation}
\end{lemma}
По этой лемме целая функция $f:=f_{\sf A}f_{\sf B}$ порядка $\ord[f]\leq 1$ с последовательностью нулей ${\sf A}\cup {\sf B}$ конечной верхней плотности при порядке $1$ удовлетворяет неравенствам 
\begin{equation*}
\ln \bigl|f(iy)\bigr|=\ln \bigl|f_{\sf A}(iy)\bigr|+\ln \bigl|f_{\sf B}(iy)\bigr|
\overset{\eqref{fA}, \eqref{fB}}\leq  v_{\alpha}(iy)+ v_{\beta}v(iy)=v(iy)
\end{equation*}
при всех $|y|\geq 1$ и  $y\in \RR$. Отсюда по  тождеству \eqref{uvm} получаем 
\begin{equation}\label{fMl}
u(iy)+\ln \bigl|f(iy)\bigr|\leq M(iy) \quad \text{при всех $|y|\geq 1$, $y\in \RR$. }
\end{equation}
При достаточно малой постоянной $a>0$ ввиду непрерывности $M$ на $i\RR$ для целой функции $h:= af$ 
конечного порядка $\ord [h]\leq 1$  с последовательностью  нулей ${\sf A}\cup {\sf B}$  можем распространить \eqref{fMl} на все $y\in \RR$ как в \eqref{vhM} с пустым $E=\varnothing$, как  в заключительном дополнении основной теоремы. Осталось показать, что $\type_1\bigl[u+\ln |h|\bigr]<+\infty$.

По построению легко видеть что заряды $\alpha-n_{\sf A}$ и $\beta -n_{\sf B}$  удовлетворяют условию Линделёфа
\eqref{con:Lp+}  рода $1$. Следовательно, и разность мер 
$$\frac{1}{2\pi}\Delta v-(n_{\sf A}+n_{\sf B})= (\alpha+\beta) -(n_{\sf A}+n_{\sf B})$$
 удовлетворяет условию Линделёфа рода $1$. По построению выше $u+v$ --- субгармоническая функция конечного типа при порядке $1$ и по теореме 
 Вейерштрасса\,--\,Адамара\,--\,Линделёфа мера Рисса 
$\frac{1}{2\pi}\Delta u+\frac{1}{2\pi}\Delta v$ удовлетворяет условию Линделёфа рода $1$.  Отсюда мера $\frac{1}{2\pi}\Delta u+(n_{\sf A}+n_{\sf B})$ удовлетворяет условию Линделёфа рода $1$. Но эта мера  конечного типа при порядке $1$ является  мерой  Рисса функции $u+\ln|h|\in \sbh_*(\CC)$ порядка $\ord[u+\ln|h|]\leq 1$. Следовательно, по теореме   Вейерштрасса\,--\,Адамара\,--\,Линделёфа\,--\,Брело функция $u+\ln |h|$ конечного типа при порядке $1$. 

Часть [Hol] с  \eqref{vhM}  для $E=\varnothing$  доказана. 

\subsection{Случай функции $M$ без ограничений вблизи мнимой оси}

Воспользуемся представлением $M:=M_{\RR}+M_{i\RR}$ из теоремы B\ref{BiR} в виде суммы двух субгармонических функций конечного типа при порядке $1$ со свойствами  \eqref{nubalR}, \eqref{nubaliR} и \eqref{vbalM} вместе с   
$\upmu \overset{\eqref{lJMmul}}{\curlyeqprec}_{|\Re|}\upmu_{\RR}$. По условию  $\upnu \curlyeqprec_{|\Re|}\upmu$. Следовательно, 
$\upnu \curlyeqprec_{|\Re|}\upmu_{\RR}$ и по \eqref{nubalR} мера Рисса субгармоническая функция  $M_{\RR}$ конечного типа при порядке $1$ гармонична вне вещественной оси. По доказанной в предыдущем подразделе части [sbh] основной теоремы для функции $M_{\RR}$  найдётся 
функция $u_{\RR}\in \sbh_*(\CC)$ конечного типа при порядке $1$ с мерой Рисса $\frac{1}{2\pi}\Delta u_{\RR}\geq \upnu$, для которой  
$u_{\RR}(iy)= M_{\RR}(iy)$ для каждого $y\in \RR$.
Отсюда для субгармонической функции  $U:=u_{\RR}+M_{i\RR}$ конечного порядка получаем 
\begin{multline}\label{MMRR}
U(iy)=u_{\RR}(iy)+M_{i\RR}(iy)\\=M_{\RR}(iy)+M_{i\RR}(iy)\overset{\eqref{{vbalM}i}}{=}M(iy)\text{ для каждого  $y\in \RR$.}
\end{multline}
Таким образом, часть [sbh]  основной теоремы доказана для любой функции $M$ конечного типа  $\type_1[M]<+\infty$ при порядке $1$. 
 
По части [Hol]  основной теоремы, применённой к функции $M_{\RR}$
из \eqref{MMRR},   
для произвольной функции $v$ с мерой Рисса $\frac{1}{2\pi}\Delta v=\upnu$
найдётся целая функция $h_{\RR}\in \Hol_*(\CC)$, для которой   функция $v+\ln |h_{\RR}|\in \sbh_*(\CC)$ конечного типа при порядке $1$ удовлетворяет ограничениям   
\begin{equation}\label{uRRMh}
v(iy)+\ln \bigl|h_{\RR}(iy)\bigr|\leq M_{\RR}(iy)\quad \text{для каждого $y\in \RR$}.
\end{equation}
Мера Рисса $\upmu_{i\RR}$ из  \eqref{nubaliR} субгармонической функции  $M_{i\RR}$ конечного типа $\type_1[M_{i\RR}]<+\infty$ сосредоточена на мнимой оси. Функция распределения $n_{\sf S}^{i\RR}:=
\lfloor\upmu_{i\RR}\rfloor$ на $\RR$, однозначно определяют последовательность   точек ${\sf S}$ на  мнимой оси со считающей целочисленной  мерой  $n_{\sf S}$ конечного типа при порядке $1$, удовлетворяющая по построению условию Линделёфа рода $1$. С помощью  \cite[леммы 2.4.3, 2.4.5]{KhDD92}
на основе стандартной техники  устанавливается  
(гораздо более общие результаты см. в \cite[следствие 2]{KhaBai16} и  \cite[теорема 2]{Kha20a}), что  
найдётся целая функция $h_{i\RR}\in \Hol_*(\CC)$
конечного типа при порядке $1$ с подпоследовательностью нулей 
${\sf S}$, для которой $\ln \bigl|h_{i\RR}(iy)\bigr|\leq M_{i\RR} (iy)$
для каждого $y\in \RR\setminus E$,
где  $E\in \mathcal B(\RR)$ удовлетворяет условию \eqref{PE}. Отсюда, если положить $h:=h_{\RR}h_{i\RR}$,  то  согласно \eqref{uRRMh} и \eqref{MMRR} получаем требуемое  \eqref{ME}.

\end{document}